\documentclass{amsart}
\usepackage{amssymb, latexsym}

\def\bi{\begin{itemize}}
\def\bs{\begin{split}}
\def\es{\end{split}}
\def\ba{\begin{align}}
\def\bas{\begin{align*}}
\def\ea{\end{align}}
\def\eas{\end{align*}}
\def\Im{{\operatorname{Im}}}

\def\Re{{\operatorname{Re}}}
\def\lo{{\operatorname{lo}}}
\def\hi{{\operatorname{hi}}}
\def\loc{{\operatorname{loc}}}

\def\C{{\mathbb C}}
\def\R{{{\mathbb R}}}
\def\Z{{{\mathbb Z}}}
\def\N{{{\mathcal N}}}
\def\S{{{\mathcal S}}}

\def\emph#1{{\it #1}}
\def\textbf#1{{\bf #1}}

\newcommand{\jx}{\langle x\rangle}
\newcommand{\eps}{\varepsilon}
\newcommand{\epso}{1+\varepsilon}
\newcommand{\epst}{3+\varepsilon}
\newcommand{\epsmo}{1-\varepsilon}
\newcommand{\intr}{\int_{\mathbb R^n}}

\theoremstyle{plain}
\newtheorem{theorem}{Theorem}
\newtheorem{definition}[theorem]{Definition}
\newtheorem{proposition}[theorem]{Proposition}
\newtheorem{lemma}[theorem]{Lemma}
\newtheorem{corollary}[theorem]{Corollary}
\newtheorem{conjecture}[theorem]{Conjecture}

\theoremstyle{remark}
\newtheorem*{remark}{Remark}

\newtheorem*{examples}{Examples}

\numberwithin{equation}{section} \numberwithin{theorem}{section}

\begin{document}

\title[Global well-posedness and scattering for NLS]
{Global well-posedness and scattering for the mass-critical nonlinear Schr\"odinger equation for radial data in high dimensions}
\author{Terence Tao}
\address{University of California, Los Angeles}
\author{Monica Visan}
\address{Institute for Advanced Study}
\author{Xiaoyi Zhang}
\address{Academy of Mathematics and System Sciences, Chinese Academy of Sciences}
\subjclass[2000]{35Q55}

\vspace{-0.3in}
\begin{abstract}
We establish global well-posedness and scattering for solutions to the defocusing mass-critical
(pseudoconformal) nonlinear Schr\"odinger equation $iu_t + \Delta u = |u|^{4/n} u$ for large
spherically symmetric $L^2_x(\R^n)$ initial data in dimensions $n\geq 3$.  After using the reductions in \cite{compact} to
reduce to eliminating blowup solutions which are almost periodic modulo scaling, we obtain a
frequency-localized Morawetz estimate and exclude a mass evacuation scenario (somewhat analogously
to \cite{ckstt:gwp}, \cite{RV}, \cite{thesis:art}) in order to conclude the argument.
\end{abstract}

\maketitle

\section{Introduction}

\subsection{The mass-critical nonlinear Schr\"odinger equation}

Fix a dimension $n \geq 1$.  We shall consider strong $L^2_x(\R^n)$ solutions to the mass-critical
(or \emph{pseudoconformal}) defocusing non-linear Schr\"odinger (NLS) equation
\begin{equation}\label{nls}
i u_t + \Delta u = F(u)
\end{equation}
where $F(u) := +|u|^{4/n} u$ is the defocusing mass-critical nonlinearity. More precisely, we say that a function $u: I \times
\R^n \to \C$ on a time interval $I \subset \R$ (possibly half-infinite or infinite) is a \emph{strong $L^2_x(\R^n)$ solution}
(or \emph{solution} for short) to \eqref{nls} if it lies in the class $C^0_{t,\loc} L^2_x(I \times \R^n) \cap L^{2(n+2)/n}_{t,\loc}
L^{2(n+2)/n}_x(I \times \R^n)$, and we have the Duhamel formula
$$
u(t_1) = e^{i(t_1-t_0)\Delta} u(t_0) - i \int_{t_0}^{t_1} e^{i(t_1-t)\Delta} F(u(t))\ dt
$$
for all $t_0, t_1 \in I$.   Here, $e^{it\Delta}$ is the propagator for the free Schr\"odinger equation defined via the Fourier transform
$$ \hat f(\xi) := \int_{\R^n} e^{-ix \cdot \xi} f(x)\ dx$$
by
$$ \widehat{e^{it\Delta} f}(\xi) = e^{-it|\xi|^2} \hat f(\xi).$$
The condition $u \in L^{2(n+2)/n}_{t,\loc} L^{2(n+2)/n}_x$ is a natural one arising from the
Strichartz perturbation theory; for instance, it is currently necessary in order to ensure
uniqueness of solutions.  Solutions to \eqref{nls} in this class have been intensively studied, see
e.g. \cite{tsutsumi}, \cite{cwI}, \cite{caz}, \cite{bourg.2d}, \cite{mt}, \cite{mv},  \cite{ck},
\cite{keraani}, \cite{begout}, \cite{tao-lens}, \cite{compact}.

It is known (see e.g. \cite{caz}) that solutions to \eqref{nls} have a conserved \emph{mass}
$$
M(u) = M(u(t)) := \int_{\R^n} |u(t,x)|^2\ dx.
$$
In particular, since our solutions lie in $L^2_x$ by definition, we have
\begin{equation}\label{massconserv}
\|u\|_{L^\infty_t L^2_x(I \times \R^n)} < \infty,
\end{equation}
whenever $u: I \times \R^n \to \C$ is a solution to \eqref{nls}.

There is a natural scaling associated to the initial value problem \eqref{nls}.  More precisely, the map
\begin{equation}\label{scaling}
\quad u(t,x)\mapsto u^\lambda(t,x):=\lambda^{-\frac n2}u\Bigl(\frac t{\lambda^2},\frac x{\lambda}\Bigr)
\end{equation}
maps a solution to \eqref{nls} to another solution to \eqref{nls} with initial data $u^\lambda(0)=u_0\bigl(\frac {x}{\lambda}\bigr)$.
The reason why \eqref{nls} is called $L^2_x$-critical (or mass-critical) is because the scaling \eqref{scaling} leaves the mass invariant.

It is known (see e.g. \cite{caz}) that if the initial data $u_0 \in L^2_x(\R^n)$ has sufficiently
small mass, then there exists a unique global solution to \eqref{nls}, which furthermore has finite
$L^{2(n+2)/n}_{t,x}(\R \times \R^n)$ norm.  This in turn implies (for either choice of sign $\pm$)
that the solution scatters to a free solution $e^{it\Delta} u_\pm$ as $t \to \pm\infty$ for some
$u_\pm \in L^2_x(\R^n)$, in the sense that
\begin{equation}\label{utm}
 \lim_{t \to \pm\infty} \| u(t)-e^{it\Delta} u_\pm \|_{L^2_x(\R^n)} = 0.
\end{equation}
Conversely, given any $u_\pm$ of sufficiently small mass there exists a solution $u$ which scatters
to it in the sense above, thus giving rise to well-defined wave and scattering operators.  See \cite{caz}
for details.

The above results were obtained by a perturbative argument and also hold in the focusing case when
the $+$ sign on the right-hand side of \eqref{nls} is replaced by a $-$ sign.  However, in the
focusing case it has long been known that large mass solutions can blow up in finite time.
Nevertheless, in the defocusing case no blowup solutions are known.  Indeed, one has the following
conjecture:

\begin{conjecture}[Global existence and scattering]\label{conj}
Given any finite mass initial data $u_0 \in L^2_x(\R^n)$, there exists a unique global solution $u
\in C^0_t L^2_x(\R \times \R^n) \cap L^{2(n+2)/n}_{t,x}(\R \times \R^n)$ to \eqref{nls} with $u(0)
= u_0$.  Furthermore, there exist $u_\pm \in L^2_x(\R^n)$ such that \eqref{utm} holds and the maps
$u_0 \mapsto u_\pm$ are homeomorphisms on $L^2_x(\R^n)$.
\end{conjecture}

\subsection{Main result}

The main result of this paper is to verify a special case of the above conjecture.

\begin{theorem}\label{main} Conjecture \ref{conj} is true when $n \geq 3$ and $u_0$ (and $u_\pm$) are restricted to be spherically symmetric.
\end{theorem}

The proof of this mass-critical theorem follows a broadly similar strategy used to settle the energy-critical
problem (see \cite{borg:scatter}, \cite{ckstt:gwp}, \cite{RV}, \cite{Monica:thesis}).  First, one
reduces to a minimal-mass blowup solution which has good localization properties in space and frequency,
establishes an initial Morawetz inequality on a frequency component, and then uses a non-critical
conservation law to prevent the solution escaping to high or low frequencies.  However, our
arguments are ``upside-down'' in the sense that the roles of high and low frequencies are reversed
from those in the energy-critical theory.  This is because the Morawetz inequality is now
subcritical instead of supercritical, and to prevent evacuation of mass to low frequencies we use the
conservation of the subcritical energy (in contrast to \cite{ckstt:gwp}, \cite{RV},
\cite{Monica:thesis}, where the conservation of the supercritical mass is used to prevent evacuation of
energy to high frequencies).

The arguments in \cite{borg:scatter}, \cite{ckstt:gwp}, \cite{RV}, \cite{Monica:thesis} were quite
quantitative, avoiding use of qualitative tools such as concentration-compactness theorems.  In this paper, we shall adopt
some qualitative technology to simplify somewhat\footnote{Roughly speaking, whereas
the quantitative approach requires managing numerous small parameters $\eta_0$, $\eta_1, \ldots$,
the qualitative approach only requires managing at most two such parameters at a time. Furthermore,
by applying limiting arguments one can often send one of the two parameters to zero or infinity.} the computations.
One consequence of this simplification is that the ``mass evacuation step'' becomes easier to
prove, as one can gain enough regularity to use the classical energy conservation law rather than a
frequency-localized variant.  The same trick can retrospectively be applied to simplify the
energy-critical theory in \cite{ckstt:gwp}, \cite{RV}, \cite{Monica:thesis}.

Our arguments rely heavily both on the high dimension $n \geq 3$ and on the spherical symmetry. The
high dimension is needed in order to enable the Morawetz inequality to have a consistent sign, and
also to make the Strichartz numerology work correctly.  The spherical symmetry is required to
localize the solution to the spatial origin (in order to be able to exploit \emph{one}-particle Morawetz
inequalities), but is also needed in order to use several powerful strengthenings of the
classical Strichartz and Sobolev inequalities, most notably the weighted Strichartz estimate of
Vilela, \cite{vilela}.  It is a challenging problem to either lower the dimension or remove the
spherical symmetry; another problem of interest would be to attack the focusing case, under the
natural additional assumption that the mass of the solution is strictly less than that of the
ground state (see \cite{merlekenig} for some recent progress on this focusing problem in the
energy-critical setting).

\textbf{Acknowledgements}: This research was partially conducted during the period Monica Visan was employed by the Clay
Mathematics Institute as a Liftoff Fellow.  This material is based upon work supported by the National Science Foundation under agreement
No. DMS--0111298.  Any opinions, findings and conclusions or recommendations expressed in this material are those of the authors and do not
reflect the views of the National Science Foundation.  The third author was supported by the NSF grant No. 10601060 (China).
The authors thank Frank Merle for helpful comments.

\section{Notation and basic estimates}\label{notation-sec}

Throughout this paper we fix the dimension $n \geq 3$.  We also fix a small exponent $\eps > 0$
depending only on $n$; for sake of concreteness, let us conservatively take $\eps :=
\frac{1}{n^{10}}$. We allow all implied constants to depend on $n$ and $\eps$.  For instance,
when we require some quantity to be sufficiently large or small, it is understood that the implied
threshold can depend on $n$ and $\eps$.

We use the notation $X \lesssim Y$, $Y \gtrsim X$, or $X = O(Y)$ to denote the estimate $|X| \leq C
Y$ for some constant $0 < C < \infty$ (which, as mentioned earlier, can depend on $n$ and $\eps$).
In some cases we shall allow the implied constant $C$ to depend on other parameters and shall
denote this by subscripts; thus, for instance, $X \lesssim_k Y$ or $X = O_k(Y)$ denotes the estimate
$|X| \leq C_k Y$ for some $C_k$ depending on $k,n,\eps$.

We use the Fourier transform to define the fractional differentiation operators $|\nabla|^s$ by the formula
$$ \widehat{|\nabla|^s f}(\xi) := |\xi|^s \hat f(\xi).$$

We shall need the following Littlewood-Paley projection operators.
Let $\varphi(\xi)$ be a bump function adapted to the ball $\{ \xi \in
\R^n: |\xi| \leq 2 \}$ which equals 1 on the ball $\{ \xi \in \R^n:
|\xi| \leq 1 \}$.  Define a \emph{dyadic number} to be any number $N \in 2^\Z$ of the form
$N = 2^j$ where $j \in \Z$ is an integer.  For each dyadic number $N$,
we define the Fourier multipliers
\begin{align*}
\widehat{P_{\leq N} f}(\xi) &:= \varphi(\xi/N) \hat f(\xi)\\
\widehat{P_{> N} f}(\xi) &:= (1 - \varphi(\xi/N)) \hat f(\xi)\\
\widehat{P_N f}(\xi) &:= (\varphi(\xi/N) - \varphi(2\xi/N)) \hat f(\xi).
\end{align*}
We similarly define $P_{<N}$ and $P_{\geq N}$.
We also define
$$ P_{M < \cdot \leq N} := P_{\leq N} - P_{\leq M} = \sum_{M < N' \leq N} P_{N'}$$
whenever $M < N$ are dyadic numbers.

The symbol $u$ shall always refer to a solution to the nonlinear Schr\"odinger equation
\eqref{nls}.  We shall use $u_N$ to denote the frequency piece $u_N := P_N u$ of $u$, and
similarly define $u_{\geq N} = P_{\geq N} u$, etc.

We use the ``Japanese bracket'' convention $\langle x \rangle := (1 + |x|^2)^{1/2}$.

We use $L^q_t L^r_x$ to denote the spacetime norm
$$ \| u \|_{L^q_t L^r_x(\R \times \R^n)} := \Bigl(\int_\R \Bigl(\int_{\R^n} |u(t,x)|^r\ dx\Bigr)^{q/r}\ dt\Bigr)^{1/r},$$
with the usual modifications when $q$ or $r$ are equal to infinity, or when the domain $\R \times \R^n$
is replaced by a smaller region of spacetime such as $I \times \R^n$.  When $q=r$ we abbreviate
$L^q_t L^q_x$ as $L^q_{t,x}$.

As we shall be frequency manipulating various Fourier multipliers it will be convenient to introduce the following definition.

\begin{definition}[H\"ormander-Mikhlin multiplier] A \emph{H\"ormander-Mikhlin multiplier} $T$ is any operator of the form
$$ \widehat{Tf}(\xi) := m(\xi) \hat f(\xi)$$
defined for any tempered distribution $f$ on $\R^n$, where the \emph{symbol} $m$ obeys the pointwise bounds
$$ |\nabla^k m(\xi)| \lesssim_{k} |\xi|^{-k}$$
for all $k \geq 0$.
\end{definition}

\begin{examples}
The Littlewood-Paley multipliers $P_N$, $P_{<N}$, $P_{\geq N}$ are H\"ormander-Mikhlin multipliers
uniformly in $N$, as are the multipliers $N^{-s} |\nabla|^s P_{<N}$ and $N^s |\nabla|^{-s} P_{\geq
N}$ for any $s \geq 0$.
\end{examples}

The classical H\"ormander-Mikhlin Theorem asserts that H\"ormander-Mikhlin multipliers are bounded
on $L^p(\R^n)$ for any $1 < p < \infty$.  We shall need an extension of this to power weights:

\begin{lemma}\label{calculus}
Let $T$ be a H\"ormander-Mikhlin multiplier, $1 < p < \infty$, and let $-\frac{n}{p} < \alpha < n - \frac{n}{p}$.  Then we have
$$ \| |x|^\alpha Tf \|_{L^p_x(\R^n)} \lesssim_{p,\alpha} \| |x|^{\alpha} f \|_{L^p_x(\R^n)}$$
for all $f$ for which the right-hand side is finite.
\end{lemma}

This estimate follows from\footnote{One can also essentially derive this estimate from the
unweighted one using Lemma \ref{bil} to control the non-local interactions when $|x| \ll |y|$ or
$|y| \ll |x|$; we omit the details.} the general Calder\'on-Zygmund theory of $A_p$ weights; see
\cite{stein:large}.

The need to deal with power weights arises primarily from our use of the following weighted Strichartz spaces.

\begin{definition}[Weighted Strichartz norms]\label{sdef}
Let $I$ be an interval and let $u: I \times \R^n \to \C$ and $G: I \times \R^n \to \C$ be
functions.  We define
$$ \| u \|_{\S(I \times \R^n)} := \| |x|^{-(1+\eps)/2} |\nabla|^{(1-\eps)/2} u \|_{L^2_{t,x}(I \times \R^n)}
    + \| u \|_{L^\infty_t L^2_x(I \times \R^n)}$$
and
$$ \| G \|_{\N(I \times \R^n)} := \| |x|^{(1+\eps)/2} |\nabla|^{-(1-\eps)/2} G \|_{L^2_{t,x}(I \times \R^n)}.$$
When the domain $I \times \R^n$ is clear from context we shall abbreviate these norms as $\S$ and $\N$, respectively.
\end{definition}

From Lemma \ref{calculus} we see that H\"ormander-Mikhlin multipliers preserve the spaces $\S$ and
$\N$.  Thus, for instance, the Littlewood-Paley multipliers are all bounded on these spaces, and one
has estimates such as
$$ \| |\nabla|^s P_{\leq N} u \|_\S \lesssim_s N^s \| u \|_\S$$
for $s \geq 0$.  Also, since the fractional integral operator $|\nabla|^{-(1-\eps)/2}$ has positive kernel, we have the comparison principle
\begin{equation}\label{comparison}
\hbox{If } G = O(H), \hbox{ then } \| G \|_{\N(I \times \R^n)} = O( \| H \|_{\N(I \times \R^n)}).
\end{equation}
We shall use the above observations in the sequel without further comment.

The relevance of these spaces to the Schr\"odinger equation arises from the following weighted Strichartz estimate of Vilela, \cite{vilela}:

\begin{proposition}[Weighted Strichartz estimates, \cite{vilela}]\label{strich}
Suppose that $u: I \times \R^n \to \C$ and $G: I \times \R^n \to \C$ solve the inhomogeneous Schr\"odinger equation $iu_t +
\Delta u = G$ in the sense of distributions.  Then
$$ \| u \|_{\S(I \times \R^n)} \lesssim \| u(t_0) \|_{L^2_x(\R^n)}
+ \| G \|_{\N(I \times \R^n)}$$
for all $t_0 \in I$.
\end{proposition}

In the spherically symmetric case the $\S$ and $\N$ norms are also related to the more traditional unweighted counterparts:

\begin{proposition}[Radial Sobolev embeddings, \cite{vilela}, \cite{stein:weiss}]\label{radsob}
If $u: I \times \R^n \to \C$ and $G: I \times \R^n \to \C$ are spherically symmetric, then
$$ \| u \|_{L^2_t L^{2n/(n-2)}_x(I \times \R^n)} \lesssim \| u \|_{\S(I \times \R^n)}$$
and
$$ \| G \|_{\N(I \times \R^n)} \lesssim \| G \|_{L^2_t L^{2n/(n+2)}_x(I \times \R^n)}.$$
\end{proposition}

\begin{proof} This follows immediately from Corollary \ref{radial embedding}.
\end{proof}

Interpolating between $L^2_t L^{2n/(n-2)}_x$ and $L^\infty_t L^2_x$ we also conclude that
\begin{equation}\label{utn}
 \| u \|_{L^{2(n+2)/n}_{t,x}(I \times \R^n)} \lesssim \| u \|_{\S(I \times \R^n)}.
\end{equation}

As a consequence of Proposition \ref{radsob}
and H\"older's inequality we immediately establish

\begin{corollary}[Basic nonlinear estimate]\label{basic}  If $u, v: I \times \R^n \to \C$ are spherically symmetric, then
\begin{align*}
\| O(|u|^{4/n} |v|) \|_{\N(I \times \R^n)}
&\lesssim  \| |u|^{4/n} |v| \|_{L^2_t L^{2n/(n+2)}_x(I \times \R^n)}  \\
&\lesssim \| u \|_{L^\infty_t L^2_x(I \times \R^n)}^{4/n} \| v \|_{L^2_t L^{2n/(n-2)}(I \times \R^n)} \\
&\lesssim \| u \|_{L^\infty_t L^2_x(I \times \R^n)}^{4/n} \| v \|_{\S(I \times \R^n)}.
\end{align*}
\end{corollary}

Standard Strichartz theory tells us that if $u$ is a solution to \eqref{nls}, then $u$ lies in
$L^2_t L^{2n/(n-2)}_x$ locally in time.  Applying Corollary \ref{basic} we see that $F(u) =
|u|^{4/n} u$ lies in $\N(I \times \R^n)$ locally in time.  Applying Proposition \ref{strich} we
thus conclude

\begin{corollary}[Local finiteness of norms]\label{local-finite}
Let $u: I \times \R^n \to \C$ be a solution to \eqref{nls}.  Then $\| u \|_{\S(J \times \R^n)} < \infty$ for all compact $J \subset I$.
\end{corollary}

This corollary will allow us to rigorously set up some continuity arguments in the sequel.

It will be important to improve upon Corollary \ref{basic} when $u$ and $v$ are separated in
frequency.  This will be accomplished by the following variant of Corollary \ref{basic}.

\begin{proposition}[Refined nonlinear estimates]\label{refined}
If $u, v: I \times \R^n \to \C$ are spherically symmetric, then we have
$$ \| |\nabla|^{\frac{1-\eps}2} O( |u|^{4/n} |v| ) \|_{\N(I \times \R^n)} \lesssim
\| |\nabla|^{\frac{n}{4}(1-\eps)} u \|_{L^\infty_t L^2_x(I \times \R^n)}^{4/n} \| |\nabla|^{-\frac{1-\eps}2} v \|_{\S(I \times \R^n)}$$
and
$$ \| |\nabla|^{\frac{1-\eps}2} O( |u|^{4/n} |v| ) \|_{\N(I \times \R^n)} \lesssim
\| |\nabla|^{\frac{3}{4}(1-\eps)} u \|_{L^\infty_t L^2_x(I \times \R^n)}^{4/n} \| |\nabla|^{(1-\eps)(\frac{1}{2} - \frac{3}{n})} v \|_{\S(I \times \R^n)}.$$
\end{proposition}

\begin{proof}  For the rest of the proof, all spacetime norms will be taken on $I \times \R^n$.
Applying Definition \ref{sdef}, the left-hand side in both inequalities can be estimated by
\begin{equation}\label{xuv}
 \| |x|^{(1+\eps)/2} |u|^{4/n} |v| \|_{L^2_{t,x}}.
\end{equation}
To prove the first estimate, we apply H\"older to bound this by
$$ \| |x|^{\frac{n}{4}(1+\eps)} u \|_{L^\infty_t L^\infty_x}^{4/n}
\| |x|^{-(1+\eps)/2} v \|_{L^2_{t,x}}$$ and the claim then follows from Definition \ref{sdef} and
Corollary \ref{radial embedding}.  To prove the second estimate, we apply H\"older slightly
differently to bound \eqref{xuv} by
$$ \| |x|^{\frac{n}{4}(1+\eps)} u \|_{L^\infty_t L^p_x}^{4/n} \| |x|^{-(1+\eps)/2} v \|_{L^2_t L^q_x}$$
where $p,q,\alpha,\beta$ are obtained by solving the equations
\begin{align*}
\frac{1}{2} - \frac{1}{q} &= \frac{1}{n} (1-\eps) (1 - \frac{3}{n}) \\
\frac{4}{np} + \frac{1}{q} &= \frac{1}{2}.
\end{align*}
One can verify that $2 \leq p,q \leq \infty$.  From Corollary \ref{radial embedding} we have
$$ \| |x|^{\frac{n}{4}(1+\eps)} u \|_{L^\infty_t L^p_x} \lesssim \| |\nabla|^{\frac{3}{4}(1-\eps)} u \|_{L^\infty_t L^2_x}$$
and
$$ \| |x|^{-(1+\eps)/2} v \|_{L^2_t L^q_x}
\lesssim \| |x|^{-(1+\eps)/2} |\nabla|^{(1-\eps)(1-\frac{3}{n})} v \|_{L^2_t L^2_x}.$$
The claim follows from Definition \ref{sdef}.
\end{proof}

\begin{remark}
Clearly, there are several more inequalities of this type; however, the above estimates are the only
ones we shall record explicitly here.
\end{remark}

\section{Overview of proof}

Let us now give an overview of the proof of Theorem \ref{main}.  By standard local well-posedness
theory (see e.g. \cite{caz}), it will suffice to prove the following quantitative estimate:

\begin{theorem}\label{main2}
Let $n \geq 3$ and let $u: I \times \R^n \to \C$ be a spherically symmetric solution to
\eqref{nls} on some time interval $I$ with the mass bound $M(u) \leq m$ for some $m < \infty$.
Then we have the spacetime bound
$$ \int_I \int_{\R^n} |u(t,x)|^{2(n+2)/n}\ dx dt \leq A(m)$$
for some finite quantity $A(m)$ depending only on $m$ (and on the dimension $n$).
\end{theorem}

\begin{remark}
In fact, Theorem \ref{main} and Theorem \ref{main2} are equivalent; see \cite{begout}, \cite{tao-lens}. We will however not need this equivalence here.
\end{remark}

We shall prove Theorem \ref{main2} by contradiction.  First we show that if Theorem~\ref{main2}
failed, then a special type of blowup solution must exist.

\begin{definition}[Almost periodic modulo scaling]\label{apdef} A solution $u: I \times \R^n \to \C$ is said to be \emph{almost periodic modulo scaling}
if there exists a function $N: I \to \R^+$ and a function $C: \R^+ \to \R^+$ such that
$$ \int_{|x| \geq C(\eta)/N(t)} |u(t,x)|^2\ dx \leq \eta$$
and
$$ \int_{|\xi| \geq C(\eta) N(t)} |\hat u(t,\xi)|^2\ d\xi \leq \eta$$
for all $t \in I$ and $\eta > 0$.
\end{definition}

\begin{remark}
The quantity $N(t)$ measures the frequency scale of the solution at time $t$.  If $u$ is not
identically zero, then $N(t)$ is uniquely defined up to a bounded multiplicative (time-dependent)
factor.  One can equivalently define $u$ to be almost periodic modulo scaling if the orbit $\{
u(t): t \in I \}$ becomes precompact in $L^2_x(\R^n)$ after quotienting out the action of the
scaling group $f(x) \mapsto \frac{1}{\lambda^{n/2}} f(\frac{x}{\lambda})$; see \cite{compact} for
further discussion.  This concept is adapted to the spherically symmetric case.  Without spherical
symmetry one also needs to take into account the translation and Galilean invariances of
\eqref{nls}, which introduce two additional modulation parameters $x(t)$ and $\xi(t)$; see
\cite{compact} for further discussion.
\end{remark}

\begin{theorem}\label{comp}
Suppose that $n \geq 1$ is such that Theorem \ref{main2} failed.  Then there exists a spherically
symmetric solution $u: I \times \R^n \to \C$, which is almost periodic modulo scaling and is such
that $\|u\|_{L^{2(n+2)/n}_{t,x}(I \times \R^n)} = +\infty$.  Furthermore, we have the frequency
bound
\begin{equation}\label{freqbound}
\sup_{t \in I} N(t) < +\infty
\end{equation}
(thus the solution does not escape to arbitrarily high frequencies).
\end{theorem}

\begin{remark}
This result will be proven in Section \ref{compact-sec}.  The result follows almost immediately
from \cite[Theorem 7.2]{compact}, but we will need an additional limiting argument in order to
extract the frequency bound \eqref{freqbound}, which we need for our argument.  This particular
component of the argument works even in low dimensions $n=1,2$, but unfortunately the remainder of
the argument relies heavily on the dimension being at least three.  Results similar to Theorem~\ref{comp}
were obtained for the energy-critical NLS in \cite{merlekenig} and for the mass-critical gKdV in \cite{merle}.
\end{remark}

In view of Theorem \ref{comp}, we see that in order to prove Theorem \ref{main2} it suffices to
show that every solution to \eqref{nls} which is almost periodic modulo scaling and obeys
\eqref{freqbound} necessarily has finite $L^{2(n+2)/n}_{t,x}$ norm.  We shall achieve this via two
key propositions in the spirit of \cite{ckstt:gwp}, \cite{RV}, \cite{Monica:thesis}.  The first
proposition establishes a frequency-localized Morawetz estimate for almost periodic solutions:

\begin{proposition}[Frequency-localized Morawetz estimate]\label{mor}
Let $n \geq 3$ and let $u: I \times \R^n \to \C$ be a spherically symmetric solution to
\eqref{nls} which is almost periodic modulo scaling and obeys \eqref{freqbound}.  Then we have
\begin{equation}\label{limn}
 \lim_{N \to \infty}
\int_I \int_{\R^n}
\frac{|\nabla u_{<N}(t,x)|^2}{|N x|^{1+\eps}}\ dx dt
 = 0.
\end{equation}
\end{proposition}

\begin{remark}
In contrast to the frequency-localized Morawetz estimates in \cite{ckstt:gwp}, \cite{RV},
\cite{Monica:thesis}, the Morawetz inequality here is of classical or ``one-particle'' type rather
than an interaction or ``two-particle'' type.  We are able to rely on one-particle inequalities due
to our assumption of spherical symmetry (cf. \cite{borg:scatter}).  Also, observe that the Morawetz
estimate here establishes spacetime control only on low-frequency components of $u$ (for any fixed
$N$), in contrast to the situation in \cite{ckstt:gwp}, \cite{RV}, \cite{Monica:thesis} in which
high-frequency components are controlled.  This is ultimately because Morawetz inequalities are
derived from variants of the momentum, which is supercritical for the energy-critical NLS, but
subcritical for the mass-critical NLS.
\end{remark}

\begin{remark}
The Morawetz inequality we use relies ultimately on the fact that $\Delta\Delta \langle x \rangle$
is non-positive, and so the argument breaks down in one and two dimensions; moreover, several other key
tools, such as the harmonic analysis estimates in Appendix \ref{appendix}, also break down in these
dimensions.  On the other hand, virial identities are equally valid in all dimensions, so it may be
that one can extend the arguments here to lower dimensions by replacing the Morawetz argument with
a virial one.
\end{remark}

\begin{remark}
With respect to the scaling \eqref{scaling}, the left-hand side of \eqref{limn} is dimensionless.
Thus, one can view this proposition as a decay estimate on the high frequencies of $u$; this is of
course consistent with the hypothesis \eqref{freqbound}.
\end{remark}

The proof of Proposition \ref{mor} is somewhat involved and will occupy Sections \ref{mor1-sec}-\ref{mor3-sec}.
The arguments used to prove Proposition~\ref{mor} automatically imply some estimates on components of $u$ in the natural solution norm
$\S(I \times \R^n)$ (see subsection~\ref{decaysec}):

\begin{proposition}[High-frequency decay of $\S$]\label{decay}
Let $n \geq 3$ and let $u: I \times \R^n \to \C$ be a spherically symmetric solution to
\eqref{nls} which is almost periodic modulo scaling and obeys \eqref{freqbound}.  Then
\begin{align}\label{finite}
\|u_{\geq N} \|_{\S(I \times \R^n)} + \frac{1}{N^{(1+\eps)/2}} \| |\nabla|^{-(1-\eps)/2} \nabla u_{< N} \|_{\S(I \times \R^n)} <\infty, \quad \forall \, N>0
\end{align}
and
\begin{align}\label{eq decay}
\lim_{N \to \infty} \Bigl[\|u_{\geq N} \|_{\S(I \times \R^n)} + \frac{1}{N^{(1+\eps)/2}} \| |\nabla|^{-(1-\eps)/2} \nabla u_{< N} \|_{\S(I \times \R^n)}\Bigr]= 0.
\end{align}
\end{proposition}

\begin{remark}
Note that \eqref{finite} follows from \eqref{eq decay}. Indeed, by \eqref{eq decay} there exists a dyadic number $N_0$ such that for $N\geq N_0$ we have
\begin{align*}
\|u_{\geq N} \|_{\S(I \times \R^n)} + \frac{1}{N^{(1+\eps)/2}} \| |\nabla|^{-(1-\eps)/2} \nabla u_{< N} \|_{\S(I \times \R^n)} \leq 1.
\end{align*}
In particular, this implies that for $N\leq N_0$
\begin{align*}
\frac{1}{N^{\frac{1+\eps}2}} \| |\nabla|^{-(1-\eps)/2} \nabla u_{< N} \|_{\S(I \times \R^n)}
&\lesssim \Bigr(\frac{N_0}N\Bigl)^{\frac{1+\eps}2} \frac{1}{N_0^{\frac{1+\eps}2}} \| |\nabla|^{-(1-\eps)/2} \nabla u_{< N_0} \|_{\S(I \times \R^n)}\\
&\lesssim \Bigr(\frac{N_0}N\Bigl)^{\frac{1+\eps}2}<\infty,
\end{align*}
which in turn implies that for $N<N_0$ we have
\begin{align*}
\| u_{N} \|_{\S(I \times \R^n)}
&\sim \frac{1}{N^{\frac{1+\eps}2}} \| |\nabla|^{-(1-\eps)/2} \nabla u_{N} \|_{\S(I \times \R^n)}\\
&\lesssim \frac{1}{N^{\frac{1+\eps}2}} \| |\nabla|^{-(1-\eps)/2} \nabla u_{< 2N} \|_{\S(I \times \R^n)}\\
&\lesssim \Bigr(\frac{N_0}N\Bigl)^{\frac{1+\eps}2}.
\end{align*}
Thus, for $N<N_0$
\begin{align*}
\|u_{\geq N} \|_{\S(I \times \R^n)}
&\lesssim \|u_{\geq N_0} \|_{\S(I \times \R^n)} + \sum_{N\leq M<N_0}\|u_{M} \|_{\S(I \times \R^n)}\\
&\lesssim 1 + \sum_{N\leq M<N_0} \Bigr(\frac{N_0}M\Bigl)^{\frac{1+\eps}2}\\
&\lesssim \Bigr(\frac{N_0}N\Bigl)^{\frac{1+\eps}2}<\infty.
\end{align*}
\end{remark}

\begin{remark}
In view of \eqref{utn}, we have now established the desired $L^{2(n+2)/n}_{t,x}$ control of $u$ on
the high frequencies.  However, the low frequencies will still require a non-trivial amount of
effort to control in this norm, even with Proposition \ref{mor} in hand.
\end{remark}

By combining Proposition \ref{decay} with a regularity argument and an energy conservation argument
(which can be viewed as a mirror image of the mass conservation argument used in \cite{ckstt:gwp},
\cite{RV}, \cite{Monica:thesis}), we shall obtain:

\begin{proposition}[Non-evacuation of mass]\label{nonevac}
Let $n \geq 3$ and let $u: I \times \R^n \to \C$ be a spherically symmetric solution to
\eqref{nls} which is almost periodic modulo scaling, is not identically zero, and obeys
\eqref{freqbound}.  Then
$$ \inf_{t \in I} N(t) > 0.$$
\end{proposition}

We prove this proposition in Section \ref{mass-sec}.
Combining Proposition~\ref{nonevac} with Proposition~\ref{decay} we now obtain

\begin{corollary}[Frequency localization implies finite lifespan]\label{cors}
Let $n \geq 3$ and let $u: I \times \R^n \to \C$ be a spherically symmetric solution to
\eqref{nls} which is almost periodic modulo scaling and obeys \eqref{freqbound}.  Suppose also that
$u$ is not identically zero.  Then $I$ is bounded.
\end{corollary}

\begin{proof}
From Proposition \ref{nonevac} we see that $N(t)$ is bounded both above and below.  From Definition
\ref{apdef} we conclude that for any $\eta > 0$ there exists $C > 0$ such that
$$
\int_{|x| \geq C} |u(t,x)|^2\ dx \leq \eta \, \hbox{ for all } t \in I.
$$
Since $u$ is non-zero, we thus see (by choosing $\eta$ small enough) that there exists $C, \eta > 0$ such that
$$
\int_{|x| < C} |u(t,x)|^2\ dx \geq \eta^{\frac{n-2}n} \, \hbox{ for all } t \in I.
$$
Thus, by H\"older,
$$
\int_{\R^n}|u(t,x)|^{\frac{2n}{n-2}}\ dx \geq c\eta \, \hbox{ for all } t \in I,
$$
for some $c>0$.  On the other hand, by Bernstein's inequality, for sufficiently small $N$ we get
$$ \int_{\R^n}|u_{< N}(t,x)|^{\frac{2n}{n-2}}\ dx \leq c\eta/10 \, \hbox{ for all } t \in I.$$
By the triangle inequality we conclude that
$$ \int_{\R^n} |u_{\geq N}(t,x)|^{\frac{2n}{n-2}}\ dx \geq c\eta/10 \, \hbox{ for all } t \in I$$
and hence, by H\"older,
$$ \| u_{\geq N} \|_{L^2_t L^{2n/(n-2)}_x(I \times \R^n)} \gtrsim |I|^{1/2}.$$
But from Proposition \ref{decay} and Proposition \ref{radsob} we know that the
left-hand side is finite, and so $I$ is bounded, as claimed.
\end{proof}

If $I$ is bounded and $\int_I \int_{\R^n} |u(t,x)|^{2(n+2)/n}\ dx dt$ is infinite, then $N(t)$
must go to infinity in finite time (see the proof of \cite[Proposition 6.1]{compact}).  But this
contradicts \eqref{freqbound}. Combining this observation with Corollary \ref{cors}, we see that any
spherically symmetric solution to \eqref{nls} which is almost periodic modulo scaling and obeys
\eqref{freqbound} must have finite $L^{2(n+2)/n}_{t,x}$ norm.  Combining this with Theorem
\ref{comp}, we obtain Theorem \ref{main2} and hence Theorem \ref{main}.

It remains to verify Theorem \ref{comp}, Proposition \ref{mor}, and Proposition \ref{nonevac}.  This
is the purpose of the remaining sections of the paper.

\section{Proof of Theorem \ref{comp}}\label{compact-sec}

We now prove Theorem \ref{comp}.  Suppose that $n \geq 1$ is such that Theorem \ref{main2} failed.
In the notation of \cite{compact}, this is precisely the assertion that the spherically symmetric
critical mass $m_{0,\operatorname{rad}}$ is finite. Thus, by \cite[Theorem 7.2]{compact} there
exists a solution $v: J \times \R^n \to \C$ to \eqref{nls} which is almost periodic modulo
scaling, and which blows up in the sense that the $L^{2(n+2)/n}_{t,x}(J \times \R^n)$ norm of $v$
is infinite.  Let $N_v(t)$ be the frequency scale function associated to $v$ as in Definition
\ref{apdef}.

We are not done yet, as $v$ does not necessarily obey the frequency bound \eqref{freqbound}.
However, we can extract a solution with this property from $v$ by a rescaling and limiting argument\footnote{These arguments are quite standard in the
literature; see e.g. \cite{merlekenig} or \cite{merle} for some other recent examples.  Our argument here is particularly close to that in
\cite{merlekenig}.}, as follows:

Write $J$ as a nested union of compact intervals $J_1 \subset J_2 \subset \ldots
\subset J$.  On each compact interval $J_i$, we have $v \in C^0_t L^2_x(J_i \times \R^n)$, which
easily implies (from Definition \ref{apdef}) that $N_v(t)$ is bounded above and below on $J_i$.
Thus, we may find $t_i \in J_i$ with the property that
\begin{equation}\label{nerv}
 N_v(t) \lesssim N_v(t_i) \hbox{ for all } t \in J_i.
\end{equation}
We choose such a time $t_i$ and then define the rescaled function $u_i: I_i \times \R^n \to \C$ as
$$ u_i(t,x) := N_v(t_i)^{n/2} v(N_v(t_i)^2 t+t_i, N_v(t_i) x )$$
where $I_i := \{ t \in \R: N_v(t_i)^2 t+t_i \in J_i \}$; thus, $I_i$ is a compact interval
containing $0$.  From the scale-invariance \eqref{scaling} and the time translation invariance, we see that
each $u_i$ is a solution to \eqref{nls}; furthermore, from Definition \ref{apdef} we see
that the initial data $\{ u_i(0): i=1,2,\ldots\}$ are a precompact subset of $L^2_x(\R^n)$.  Thus,
by passing to a subsequence if necessary, we may assume that $u_i(0)$ converges strongly in
$L^2_x(\R^n)$ to another initial datum $u_0 \in L^2_x(\R^n)$.  Also, from the conservation of mass we
know that $u_i(0)$ all have the same $L^2_x(\R^n)$ norm, which is non-zero as $v$ is not identically zero.
Thus $u_0$ is also not identically zero.

Let $\tilde u_i: \tilde I_i \times \R^n \to \C$ be the maximal Cauchy extension of $u_i$; thus,
$\tilde u_i$ is the maximal-lifespan solution to \eqref{nls} which agrees with $u_i$ on $I_i$. Let
$u: (-T_-, T_+) \times \R^n \to \C$ be the maximal Cauchy development to \eqref{nls} with initial
data $u_0$ for some $-\infty \le -T_- < 0 < T_+ \le +\infty$.  If $I'$ is any compact subinterval
of $(-T_-,T_+)$ containing $0$, then $u$ has finite $L^{2(n+2)/n}_{t,x}$ norm on $I' \times \R^n$.
Since $\tilde u_i(0)$ converges to $u(0)$ strongly in $L^2_x(\R^n)$, we thus see from the standard
well-posedness theory (see e.g. \cite{compact}) that for sufficiently large $i$, $\tilde I_i$
contains $I'$, and $\tilde u_i$ converges in $C^0_t L^2_x( I' \times \R^n ) \cap
L^{2(n+2)/n}_{t,x}(I' \times \R^n)$ to $u$.  In particular, $\tilde u_i$ has bounded
$L^{2(n+2)/n}_{t,x}(I' \times \R^n)$ norm as $i \to \infty$.  On the other hand, from the monotone
convergence theorem we have
$$ \| v \|_{L^{2(n+2)/n}_{t,x}(J_i \times \R^n)} \to \infty \hbox{ as } i \to \infty,$$
which after rescaling becomes
$$ \| u_i \|_{L^{2(n+2)/n}_{t,x}(I_i \times \R^n)} \to \infty \hbox{ as } i \to \infty.$$
The only way these facts can be consistent is if $I_i \not \subseteq I'$ for all sufficiently large
$i$.  But $I'$ was an arbitrary subinterval of $(-T_-,T_+)$ containing $0$.  After passing to a
subsequence if necessary (and using the usual diagonalization trick), this leaves only two
possibilities:
\begin{itemize}
\item For every $0 < t < T_+$, $I_i$ contains $[0,t]$ for all sufficiently large $i$.
\item For every $-T_- < t < 0$, $I_i$ contains $[-t,0]$ for all sufficiently large $i$.
\end{itemize}
By time reversal symmetry, it suffices to consider the former possibility.  Then, for any $0 \leq t <
T_+$ we see that $u(t)$ can be approximated to arbitrary accuracy in the $L^2_x(\R^n)$ norm by
$u_i(t)$, which is a rescaled version of a function in the orbit $\{ v(t): t \in J \}$.  But the
latter set is precompact in $L^2_x(\R^n)$ after quotienting out by scaling.  Thus, the orbit $\{
u(t): 0 \leq t< T_+ \}$ is also precompact in $L^2_x(\R^n)$ after quotienting out by scaling.  In
other words, if we set $I := [0,T_+)$ then $u: I \times \R^n \to \C$ is almost periodic modulo
scaling.

We now claim that $u$ blows up.  For if $u$ had finite $L^{2(n+2)/n}_{t,x}(I \times \R^n)$ norm,
then (since $(-T_-,T_+)$ was the maximal Cauchy development) the standard local well-posedness
theory (see e.g. \cite{caz}) would imply that $T_+ = +\infty$ and that $u$ scattered to a free
solution $e^{it\Delta} u_+$ as $t \to +\infty$, that is, $\lim_{t \to +\infty} \| u(t) - e^{it\Delta}
u_+ \|_{L^2_x(\R^n)} = 0$.  But a stationary phase (or fundamental solution) analysis of this free
solution reveals that this scattering is only compatible with the almost periodicity of $u$ modulo
scaling if $u_+ = 0$ (cf. \cite{compact}).  Conservation of mass then forces $u$ to be identically
zero, a contradiction.  Hence, $u$ blows up.

Finally, we need to show \eqref{freqbound}.  Let $\eta > 0$ be arbitrary.  From
\eqref{nerv} and Definition~\ref{apdef} there exists $C(\eta) > 0$ such that
$$ \| P_{> C(\eta) N_v(t_i)} v(t) \|_{L^2_x(\R^n)} \leq \eta$$
for all $i$ and all $t \in J_i$.  Rescaling this, we obtain
$$ \| P_{> C(\eta) } u_i(t) \|_{L^2_x(\R^n)} \leq \eta$$
for all $i$ and all $t \in I_i$.  Since $u_i$ converges strongly in $C^0_t L^2_x$ to $u$ on $[0,t]$ for any $0 < t < T_+$, we conclude that
$$ \| P_{> C(\eta) } u(t) \|_{L^2_x(\R^n)} \leq \eta$$
for all $0 < t < T_+$.  Comparing this to Definition \ref{apdef} (and the fact that $u$ has
non-zero mass), we conclude \eqref{freqbound} as desired.  This proves Theorem \ref{comp}.

\section{Proof of the Morawetz inequality I. Scaling}\label{mor1-sec}
We now turn to the proof of Proposition \ref{mor}.  We begin by using a scaling argument to
eliminate the role of the frequency parameter $N$.  We first give a simple high-frequency mass
decay estimate that follows from \eqref{freqbound} (compare with Proposition \ref{decay}).

\begin{lemma}[Mass decay at high frequencies]\label{massdecay}
Let $u: I \times \R^n \to \C$ be a solution to \eqref{nls} which is almost periodic modulo scaling
and obeys \eqref{freqbound}.  Then,
$$ \| u \|_{L^\infty_t L^2_x(I \times \R^n)} < \infty$$
and
$$ \lim_{N \to \infty} \| u_{\geq N} \|_{L^\infty_t L^2_x(I \times \R^n)}
+ \frac{1}{N^{(1+\eps)/2}} \||\nabla|^{-(1-\eps)/2} \nabla u_{< N} \|_{L^\infty_t L^2_x(I \times \R^n)} = 0.$$
\end{lemma}

\begin{proof}
The first bound is just \eqref{massconserv}, so we turn to the second bound.  Let $\eta > 0$ be
arbitrary.  Then, from \eqref{freqbound} and Definition \ref{apdef} we see that there exists
$C(\eta) > 0$ such that
$$ \| P_{>C(\eta)} u(t) \|_{L^2_x(\R^n)} \leq \eta$$
for all $t \in I$.  This already gives
\begin{align}\label{decay mass hf}
\lim_{N \to \infty} \| u_{\geq N} \|_{L^\infty_t L^2_x(I \times \R^n)} = 0.
\end{align}
For the second term, we split $u_{<N} = u_{<\sqrt{N}} + u_{\sqrt{N} \leq \cdot < N}$ and compute
$$
\frac{1}{N^{\frac{1+\eps}2}} \||\nabla|^{-\frac{1-\eps}2} \nabla u_{< N} \|_{L^\infty_t L^2_x(I \times \R^n)}
\lesssim \frac{1}{N^{\frac{1+\eps}4}} \| u \|_{L^\infty_t L^2_x(I \times \R^n)}
+ \| u_{\ge \sqrt{N}} \|_{L^\infty_t L^2_x(I \times \R^n)}.$$
Using \eqref{decay mass hf}, we see that the right-hand side of the above inequality goes to zero as $N \to \infty$, as claimed.
\end{proof}

In view of this lemma, we see that Proposition \ref{mor} will follow from the following variant,
which does not explicitly assume almost periodicity modulo scaling.

\begin{proposition}[Frequency-localized Morawetz estimate, reformulated]\label{mor2}\leavevmode\\
Let $n \geq 3$, $m > 0$, and $0 < \eta < 1$.  Then, there exists $\delta > 0$ with the following
property: given any $N > 0$ and any spherically symmetric solution $u: I \times \R^n \to \C$ to
\eqref{nls} which obeys the bounds
$$ \| u \|_{L^\infty_t L^2_x(I \times \R^n)} \leq m$$
and
$$ \| u_{\geq N} \|_{L^\infty_t L^2_x(I \times \R^n)}
+ \frac{1}{N^{(1+\eps)/2}} \| |\nabla|^{-(1-\eps)/2}\nabla u_{< N} \|_{L^\infty_t L^2_x(I \times \R^n)} \leq \delta,$$
we have
\begin{equation}\label{qi}
\int_I \int_{\R^n} \frac{|\nabla u_{<N}(t,x)|^2}{|N x|^{1+\eps}}\ dx dt
\leq \eta.
\end{equation}
\end{proposition}

The point of reformulating Proposition \ref{mor} in this way is that the scale invariance
\eqref{scaling} does not affect any component of the hypothesis or conclusion, other than by
changing $I$ and $N$.  Thus we may normalize $N=1$. By a limiting argument, we may then take $I$ to
be compact. Now we observe that by Corollary~\ref{local-finite}, the left-hand side of \eqref{qi} varies
continuously in $I$ and goes to zero when $I$ shrinks to a point.  Thus, by standard continuity
arguments, it suffices to show the following bootstrap version of the proposition:

\begin{proposition}[Frequency-localized Morawetz estimate, normalized bootstrap version]\label{mor3}\leavevmode\\
Let $n \geq 3$, $m > 0$, and $0 < \eta < 1$.  Then, there exists $\delta > 0$ with the following
property: given any spherically symmetric solution $u: I \times \R^n \to \C$ to \eqref{nls} with
$I$ compact, which obeys the mass bound
\begin{equation}\label{mass-bound}
 \| u \|_{L^\infty_t L^2_x(I \times \R^n)} \leq m
\end{equation}
and the high-frequency decay bound
\begin{equation}\label{delta-bound}
 \| u_{\hi} \|_{L^\infty_t L^2_x(I \times \R^n)}
+ \||\nabla|^{-(1-\eps)/2} \nabla u_{\lo} \|_{L^\infty_t L^2_x(I \times \R^n)} \leq \delta,
\end{equation}
where $u_\hi := u_{\geq 1}$ and $u_\lo := u_{<1}$, such that we also have the bootstrap hypothesis
\begin{equation}\label{boot}
Q_I \leq 2 \eta,
\end{equation}
where $Q_I$ is the quantity
\begin{equation}\label{qidef}
 Q_I :=
\int_I \int_{\R^n} \frac{|\nabla u_{\lo}(t,x)|^2}{|x|^{1+\eps}}\ dx dt,
\end{equation}
then we have
$$ Q_I \leq \eta.$$
\end{proposition}

It remains to prove Proposition \ref{mor3}.  This will occupy the next two sections of this paper.

\section{Proof of the Morawetz inequality II. High and low frequency estimates}\label{mor2-sec}

In this section we exploit the hypotheses \eqref{mass-bound}, \eqref{delta-bound}, \eqref{boot} to
establish some estimates on $u_\lo$ and $u_\hi$.  We begin with the low-frequency estimates.
Throughout this section we omit the domain $I \times \R^n$ for brevity.

\begin{proposition}[Low-frequency estimates]\label{lowfreq}  Let the hypotheses be as in
Proposition \ref{mor3}.  Then,
\begin{align}
\| |\nabla|^{-(1-\eps)/2} \nabla u_\lo \|_{\S} &\lesssim \eta^{1/2} \label{lo-so}\\
\| \nabla u_\lo \|_{\S} &\lesssim \eta^{1/2} \label{lo-so2}\\
\| \nabla u_\lo \|_{L^2_t L_x^{2n/(n-2)}} &\lesssim \eta^{1/2}.
\label{lo-strich}
\end{align}
\end{proposition}

\begin{proof}  We begin with \eqref{lo-so}.  From Definition \ref{sdef}, we need to establish
$$ \| |\nabla|^{-(1-\eps)/2} \nabla u_\lo \|_{L_t^\infty L^2_{x}} \lesssim \eta^{1/2}$$
and
$$ \| |x|^{-(1+\eps)/2}  \nabla u_\lo \|_{L^2_{t,x}} \lesssim \eta^{1/2}.$$
The first claim follows from \eqref{delta-bound}, taking $\delta=\delta(\eta)$ sufficiently small.
The second claim follows from \eqref{boot} and \eqref{qidef}.

The claim \eqref{lo-so2} follows from \eqref{lo-so} by writing $\nabla u_\lo =
|\nabla|^{\frac{1-\eps}2} P_{<100} |\nabla|^{-\frac{1-\eps}2} \nabla u_\lo$ and recalling that the
H\"ormander-Mikhlin multiplier $|\nabla|^{\frac{1-\eps}2} P_{<100}$ is bounded on $\S$.

The claim \eqref{lo-strich} follows from \eqref{lo-so2} and Proposition~\ref{radsob}.
\end{proof}

\begin{remark}
Note that \eqref{lo-strich} and Corollary~\ref{radial embedding} (or Hardy's inequality) implies that in dimensions $n>4$ we have
\begin{align}\label{lo-hardy}
\| |x|^{-1} u_\lo \|_{L^2_t L_x^{2n/(n-2)}} &\lesssim \eta^{1/2}.
\end{align}
In lower dimensions, that is $n=3,4$, we cannot expect such a strong decay for the low frequencies.  However, by Corollary~\ref{radial embedding}
(or Hardy's inequality), Sobolev embedding, interpolation, Bernstein, \eqref{delta-bound}, and \eqref{lo-so}, we get the following decay estimate
which is valid in all dimensions $n\geq 3$ and sufficient for our purposes:
\begin{align}
\||x|^{-\frac{n(1-\eps)}{2(n+4)}} u_\lo\|_{L^{\frac{2(n+4)}n}_{t,x}}
&\lesssim \||\nabla|^{\frac{n(1-\eps)}{2(n+4)}} u_\lo\|_{L^{\frac{2(n+4)}n}_{t,x}} \notag\\
&\lesssim \||\nabla|^{\frac{n(3-\eps)}{2(n+4)}} u_\lo\|_{L^{\frac{2(n+4)}n}_t L^{\frac{2(n+4)}{n+2}}_x}\notag\\
&\lesssim \||\nabla|^{\frac{n(3-\eps)}{2(n+4)}} u_\lo\|_{L^\infty_t L^2_x}^{\frac{4}{n+4}} \||\nabla|^{\frac{n(3-\eps)}{2(n+4)}} u_\lo\|_{L^2_t L^{\frac{2n}{n-2}}_x}^{\frac n{n+4}} \notag\\
&\lesssim \||\nabla|^{-\frac{1-\eps}2} \nabla u_\lo\|_{L^\infty_t L^2_x}^{\frac{4}{n+4}} \||\nabla|^{-\frac{1-\eps}2} \nabla u_\lo\|_{L^2_t L^{\frac{2n}{n-2}}_x}^{\frac n{n+4}} \notag\\
&\lesssim \delta^{\frac 4{n+4}}\||\nabla|^{-\frac{1-\eps}2} \nabla u_\lo\|_{\S}^{\frac n{n+4}}\notag\\
&\lesssim \delta^{\frac 4{n+4}}.\label{lo-sd}
\end{align}
\end{remark}

Now we establish high-frequency estimates.

\begin{proposition}[High-frequency estimates]\label{highfreq}
Let the hypotheses be as in Proposition \ref{mor3}.  If $\delta$ is sufficiently small, then
$$ \| u_\hi \|_{\S} \lesssim \delta +\delta^{4/n}.$$
\end{proposition}

\begin{proof}
Applying $P_\hi := P_{\geq 1}$ to \eqref{nls}, we see that
$$ (i\partial_t + \Delta) u_\hi = P_\hi F(u).$$
From Proposition \ref{strich} and \eqref{mass-bound}, we conclude
$$ \| u_\hi \|_{\S} \lesssim \delta + \| P_\hi F(u) \|_{\N}.$$
We then split
\begin{align*}
P_\hi F(u) &= P_\hi F(u_\lo) + P_\hi( F(u) - F(u_\lo) )  \\
&= \Delta^{-1} \nabla \cdot P_\hi (\nabla F(u_\lo)) + P_\hi O( |u_\lo|^{4/n} |u_\hi| + |u_\hi|^{1+4/n} ) \\
&= \Delta^{-1} \nabla \cdot P_{>1/100} P_\hi O( |u_\lo|^{4/n} |\nabla u_\lo| )+ P_\hi O( |u_\lo|^{4/n} |u_\hi| )\\
&\quad + P_\hi O( |u_\hi|^{1+4/n} ).
\end{align*}
Discarding the H\"ormander-Mikhlin multipliers $\Delta^{-1} \nabla \cdot P_{>1/100}$ and $P_\hi$,
we thus conclude from this and \eqref{comparison} that
$$ \| u_\hi \|_{\S} \lesssim \delta + \| P_\hi( |u_\lo|^{4/n} |\nabla u_\lo| ) \|_{\N}
+
\| P_\hi( |u_\lo|^{4/n} |u_\hi| ) \|_{\N}
+
\| |u_\hi|^{4/n} |u_\hi| \|_{\N}.$$
From Corollary \ref{basic} and \eqref{delta-bound}, we have
$$
\| |u_\hi|^{4/n} |u_\hi| \|_{\N}
\lesssim \delta^{4/n} \| u_\hi \|_{\S}.$$
Meanwhile, if we discard the  H\"ormander-Mikhlin multiplier
$P_\hi |\nabla|^{-(1-\eps)/2}$ and apply the first estimate in Proposition \ref{refined}, we get
\begin{align*}
\| P_\hi( |u_\lo|^{4/n} |u_\hi| ) \|_{\N}
&\lesssim
\| |\nabla|^{(1-\eps)/2}( |u_\lo|^{4/n} |u_\hi| ) \|_{\N} \\
&\lesssim
\| |\nabla|^{\frac{n}{4}(1-\eps)} u_\lo \|_{L^\infty_t L^2_x}^{4/n} \| |\nabla|^{-(1-\eps)/2} u_\hi \|_{\S} \\
&\lesssim
\| |\nabla|^{\frac{n}{4}(1-\eps)} u_\lo \|_{L^\infty_t L^2_x}^{4/n} \| u_\hi \|_{\S}.
\end{align*}
A similar argument using \eqref{lo-so} gives
$$
\| P_\hi( |u_\lo|^{4/n} |\nabla u_\lo| ) \|_{\N}
\lesssim \eta^{1/2}\| |\nabla|^{\frac{n}{4}(1-\eps)} u_\lo \|_{L^\infty_t L^2_x}^{4/n}.$$
Putting all these together, we obtain
$$
\| u_\hi \|_{\S} \lesssim (\delta + \delta^{4/n} +
\| |\nabla|^{\frac{n}{4}(1-\eps)} u_\lo \|_{L^\infty_t L^2_x}^{4/n}) (1 + \| u_\hi \|_{\S}).$$
As from \eqref{delta-bound} and Bernstein we have
$$ \| |\nabla|^{\frac{n}{4}(1-\eps)} u_\lo \|_{L^\infty_t L^2_x}^{4/n} \lesssim \delta^{4/n},$$
we obtain
$$ \| u_\hi \|_{\S} \lesssim (\delta+\delta^{4/n}) + (\delta+\delta^{4/n}) \| u_\hi \|_{\S}.$$
From Corollary \ref{local-finite} and the compactness of $I$, the norm appearing on both sides of
this estimate is finite, so for $\delta$ small enough we obtain the claim.
\end{proof}

\section{Proof of the Morawetz inequality III. Monotonicity formula}\label{mor3-sec}

To conclude the proof of Proposition \ref{mor3} we shall need a monotonicity formula that gives a
nontrivial estimate on the spacetime integral $Q_I$.  We shall phrase this monotonicity formula in
the context of a general forced NLS, as follows.

\begin{proposition}[General Morawetz inequality]\label{genmor}  Let $I$ be an interval, let $n \geq 3$, and let
$\phi, G \in C^0_t H^1_x(I \times \R^n)$ solve the equation
\begin{equation}\label{forced}
i \phi_t + \Delta \phi = F(\phi) + G
\end{equation}
on $I \times \R^n$.  Let $\eps > 0$.  If $\eps$ is sufficiently small depending on $n$, then we have
\begin{align*}
&\int_I \int_{\R^n} \Bigl(\frac{|\phi(t,x)|^2}{\langle x \rangle^{3+\eps}}
+ \frac{|\phi(t,x)|^{2(n+2)/n}}{\langle x \rangle} + \frac{|\nabla \phi(t,x)|^2}{\langle x \rangle^{1+\eps}}\Bigr)\ dx dt \\
&\quad \lesssim_\eps \sup_{t \in I} \| \phi(t) \|_{L^2_x(\R^n)} \| \nabla \phi(t) \|_{L^2_x(\R^n)} \\
&\quad\quad + \int_I \int_{\R^n} |G(t,x)| |\nabla \phi(t,x)|\ dx
dt\\
&\quad\quad + \int_I \int_{\R^n} \frac{1}{\langle x\rangle} |G(t,x)| |\phi(t,x)|\ dx dt.
\end{align*}
\end{proposition}

\begin{remark}
This is not the sharpest Morawetz inequality we can establish; for instance, one can replace $\|
\phi(t)\|_{L^2_x(\R^n)} \| \nabla \phi(t)\|_{L^2_x(\R^n)}$ on the right-hand side by $\|\phi(t)
\|_{\dot H^{1/2}_x(\R^n)}^2$, and in dimensions $n \geq 4$ one can remove the $\eps$ in the first
denominator on the left-hand side.  One can also lower the regularity required on $\phi$ and $G$.
However, the estimate as stated is sufficient for our purposes; in fact, only the last term on the
left-hand side will actually be used.  Note that this estimate does not require spherical symmetry;
however, as the estimate is localized to the spatial origin, it is not particularly effective in
the general (translation-invariant) setting in which the assumption of spherical symmetry is
dropped.
\end{remark}

\begin{proof}
By standard limiting arguments we may assume that $I$ is a compact interval, and $\phi$ is smooth
in time and Schwartz in space.  We introduce the spatial weight
$$ a(x) := \langle x \rangle - \eps \langle x \rangle^{1-\eps}$$
and consider the Morawetz functional $M_a: I \to \R$ defined by
$$ M_a(t) := 2 \int_{\R^n}
\nabla a(x)\Im(\bar\phi(t,x)\nabla \phi(t,x))dx.
$$
Since $\nabla a = O(1)$, we see from Cauchy-Schwartz that
$$ | M_a(t)| \lesssim \| \phi(t) \|_{L^2_x(\R^n)} \| \nabla \phi(t) \|_{L^2_x(\R^n)}$$
and hence, by the Fundamental Theorem of Calculus,
$$ \int_I \partial_t M_a(t)\ dt \lesssim \sup_{t \in I} \| \phi(t) \|_{L^2_x(\R^n)} \| \nabla \phi(t) \|_{L^2_x(\R^n)}.$$
To establish the proposition it thus suffices to show that
\begin{align*}
\partial_t M_a(t) &\geq c_\eps
\int_{\R^n} \Bigl(\frac{|\phi(t,x)|^2}{\langle x \rangle^{3+\eps}}
+ \frac{|\phi(t,x)|^{2(n+2)/n}}{\langle x \rangle} + \frac{|\nabla \phi(t,x)|^2}{\langle x \rangle^{1+\eps}}\Bigr)\ dx \\
&\quad -c\int_{\R^n} \Bigl(|G(t,x)| |\nabla \phi(t,x)|- \frac{1}{\langle x\rangle} |G(t,x)| |\phi(t,x)|\Bigr)\ dx.
\end{align*}

A direct calculation establishes that
\begin{align*}
\partial_t M_a(t)&=\intr (-\Delta \Delta a(x))|\phi(t,x)|^2dx+4\intr
a_{jk}(x)\Re(\bar\phi_j(t,x)\phi_k(t,x))dx\\
&\ \ +2\intr \nabla a(x)\{F(\phi) + G,\phi\}_p(t,x)\ dx
\end{align*}
where $\{,\}_p$ denotes the \emph{momentum bracket}
$$ \{f,g\}_p  := \Re( f \overline{\nabla g} - g \overline{\nabla f} ).$$

Thus, it suffices to establish the estimates
\begin{gather}
\intr (-\Delta \Delta a(x))|\phi(t,x)|^2\ dx \gtrsim_\eps
    \int_{\R^n} \frac{|\phi(t,x)|^2}{\langle x \rangle^{3+\eps}}\ dx \label{minor1} \\
\intr a_{jk}(x)\Re(\bar\phi_j(t,x)\phi_k(t,x))\ dx \gtrsim_\eps
    \int_{\R^n} \frac{|\nabla \phi(t,x)|^2}{\langle x \rangle^{1+\eps}}\ dx \label{minor2} \\
\intr \nabla a(x)\{F(\phi),\phi\}_p(t,x)\ dx \gtrsim_\eps
\int_{\R^n} \frac{|\phi(t,x)|^{2(n+2)/n}}{\langle x \rangle}\ dx\label{minor3} \\
\Bigl|\intr\nabla a(x)\{G,\phi\}_p(t,x)\ dx\Bigr|
    \lesssim \int_{\R^n} |G(t,x)| |\nabla \phi(t,x)| + \frac{|G(t,x)| |\phi(t,x)|}{\langle x\rangle} \ dx. \label{minor4}
\end{gather}
To achieve this, we compute
\begin{align}
a_j(x)&=\frac {x_j}{\jx}-\eps(\epsmo)\frac {x_j}{\jx^{\epso}},\nonumber\\
a_{jk}(x)&=\frac {\delta_{jk}}{\jx}-\frac
{x_jx_k}{\jx^3}-\eps(\epsmo)\frac {\delta_{jk}}{\jx^{\epso}}+
\eps(\epsmo)(\epso)\frac {x_jx_k}{\jx^{\epst}},\label{aderiv}\\
\Delta a(x)&=\frac {n}{\jx}-\frac
{|x|^2}{\jx^3}-\eps(\epsmo)\frac {n}{\jx^{\epso}}+
\eps(\epsmo)(\epso)\frac {|x|^2}{\jx^{\epst}},\label{adelta}\\
-\Delta\Delta a(x)&=\frac {(n-1)(n-3)}{\jx^3}-\frac
{\eps(\epsmo)(\epso)(n-1-\eps)(n-3-\eps)}{\jx^{3+\eps}}\\
&\ \ +\frac {6(n-3)}{\jx^5}-\frac {2\eps(\epsmo)(\epso)(\epst)(n-3-\eps)}{\jx^{5+\eps}}\nonumber
\\
&\ \ +\frac {15}{\jx^7}-\frac {\eps(\epsmo)(\epso)(\epst)(5+\eps)}{\jx^{7+\eps}}. \nonumber
\end{align}
For $\eps$ sufficiently small we now see that
$$ -\Delta \Delta a(x) \gtrsim_\eps \frac{1}{\langle x \rangle^{3+\eps}},$$
which gives \eqref{minor1}.

To prove \eqref{minor2}, it suffices (by splitting $\phi$ into real and imaginary parts) to establish the pointwise estimate
$$
a_{jk}(x) v_j v_k \gtrsim \frac{|v|^2}{\langle x \rangle^{1+\eps}}
$$
for any $x \neq 0$ and any real vector $v \in \R^n$.  We expand the left-hand side using \eqref{aderiv} as
$$ |v|^2 \bigl( \langle x \rangle^{-1} - \eps(1-\eps) \langle x \rangle^{-1-\eps} \bigr) - \Bigl| \frac{x}{|x|} \cdot v\Bigr|^2
|x|^2 \bigl( \langle x \rangle^{-3} - \eps(1-\eps)(1+\eps) \langle x \rangle^{-3-\eps} \bigr).$$
Since $| \frac{x}{|x|} \cdot v|$ ranges between $0$ and $|v|$, it thus suffices to show that
$$ \langle x \rangle^{-1} - \eps(1-\eps) \langle x \rangle^{-1-\eps} \gtrsim \langle x \rangle^{-1-\eps}$$
and
$$  \langle x \rangle^{-1} - \eps(1-\eps) \langle x \rangle^{-1-\eps}
- |x|^2 \langle x \rangle^{-3} +  \eps(1-\eps)(1+\eps) |x|^2 \langle x \rangle^{-3-\eps}  \gtrsim \langle x \rangle^{-1-\eps}.$$
The first claim is clear when $\eps$ is sufficiently small.  To see the second, we use the estimates
$$  \langle x \rangle^{-1} - |x|^2 \langle x \rangle^{-3} = \langle x \rangle^{-3}$$
and
$$ -\eps(1-\eps)(1+\eps) \langle x \rangle^{-1-\eps}
+ \eps(1-\eps)(1+\eps) |x|^2 \langle x \rangle^{-3-\eps}  = O(\eps \langle x \rangle^{-3-\eps})$$
to rewrite the left-hand side as
$$ \langle x \rangle^{-3} + \eps^2 (1-\eps) \langle x \rangle^{-1-\eps}
+ O(\eps \langle x \rangle^{-3-\eps})$$
and the claim is now clear.

Next, we establish \eqref{minor3}.  Observe the identity
$$ \{ F(\phi), \phi \}_p = - \frac{2}{n+2} \nabla |\phi|^{\frac{2(n+2)}{n}}.$$
Integrating by parts, the left-hand side of \eqref{minor3} becomes
$$ \frac{2}{n+2} \intr \Delta a(x) |\phi|^{2(n+2)/n}\ dx.$$
But from \eqref{adelta} we see that $\Delta a(x) \gtrsim_\eps \langle x\rangle^{-1-\eps}$, and the claim follows.

Finally, we establish \eqref{minor4}.  Observe that
$$ \{G,\phi\}_p =  2 \Re( G \overline{\nabla \phi} ) - \nabla \Re( \phi \overline{G} ).$$
Integrating by parts and using the crude bounds $\nabla a = O(1)$,
$\nabla^2 a = O( 1/\langle x \rangle )$ we obtain the claim.
\end{proof}

\emph{Proof of Proposition~\ref{mor3}}.

By applying $P_{\lo} := P_{<1}$ to \eqref{nls}, we see that $\phi := u_\lo$ will solve \eqref{forced} with $G$ equal to the nonlinear commutator
\begin{equation}\label{commutator}
 G := P_{\lo} F(u) - F( P_\lo u ).
\end{equation}
Applying the hypotheses \eqref{mass-bound}, \eqref{delta-bound}, and Bernstein, we conclude from Proposition~\ref{genmor} that
$$
\int_I \int_{\R^n} \frac{|\nabla u_\lo (t,x)|^2}{\langle x \rangle^{1+\eps}}\ dx dt \lesssim_{m,\eps} \delta
    + \int_I \int_{\R^n} |G(t,x)| \Bigl(|\nabla u_\lo(t,x)| + \frac{|u_\lo(t,x)|}{\langle x \rangle}\Bigr) \ dx dt.$$
By the uncertainty principle (Lemma \ref{smoother}), the left-hand side controls $Q_I$. Thus, to
conclude the proof of Proposition \ref{mor3} (and hence Proposition \ref{mor}) it will suffice to
establish the bound
$$
\int_I \int_{\R^n} |G(t,x)| \Bigl(|\nabla u_\lo(t,x)| + \frac{|u_\lo(t,x)|}{\langle x \rangle}\Bigr) \ dx dt \lesssim_{m} \delta^\eps
$$
for $\delta$ sufficiently small depending on $m$ and $\eps$, since we may then take $\delta$ small compared to $\eta$.

By H\"older and \eqref{lo-strich}, we estimate
\begin{align*}
\int_I \int_{\R^n} |G(t,x)| |\nabla u_\lo(t,x)| \ dx dt
&\lesssim \|G\|_{L^2_t L^{2n/(n+2)}_x(I \times \R^n)}\|\nabla u_\lo\|_{L^2_t L^{2n/(n-2)}_x(I \times \R^n)}\\
&\lesssim \|G\|_{L^2_t L^{2n/(n+2)}_x(I \times \R^n)},
\end{align*}
while in dimensions $n>4$ by H\"older and \eqref{lo-hardy}, we estimate
\begin{align*}
\int_I \int_{\R^n} |G(t,x)| \frac{|u_\lo(t,x)|}{\langle x \rangle} \ dx dt
&\lesssim \|G\|_{L^2_t L^{2n/(n+2)}_x(I \times \R^n)}\Bigl\|\frac{u_\lo}{\langle x \rangle}\Bigr\|_{L^2_t L^{2n/(n-2)}_x(I \times \R^n)}\\
&\lesssim \|G\|_{L^2_t L^{2n/(n+2)}_x(I \times \R^n)}.
\end{align*}
Thus, in dimensions $n>4$ we reduce to showing that
$$ \| G \|_{L^2_t L^{2n/(n+2)}_x(I \times \R^n)} \lesssim_{m} \delta^\eps.$$
We split the commutator \eqref{commutator} as
\begin{align*}
G &= P_\lo [F(u) - F( u_\lo )] - P_\hi( F(u_\lo) ) \\
&= P_\lo O( |u_\hi| |u_\lo|^{4/n} + |u_\hi|^{1+4/n} ) - \Delta^{-1} \nabla \cdot P_\hi( \nabla F(u_\lo) ) \\
&= P_\lo O(|u_\hi| |u_\lo|^{4/n}) + P_\lo O( |u_\hi|^{1+4/n} )
+ \Delta^{-1} \nabla \cdot P_\hi O( |u_\lo|^{4/n} |\nabla u_\lo| ).
\end{align*}
As the multipliers $P_\lo$ and $\Delta^{-1} \nabla \cdot P_\hi$ are bounded on
$L^{2n/(n+2)}_x(\R^n)$, we reduce to showing that
\begin{align*}
\| |u_\hi| |u_\lo|^{4/n} \|_{L^2_t L^{2n/(n+2)}_x(I \times \R^n)}
&+ \| |u_\hi| |u_\hi|^{4/n} \|_{L^2_t L^{2n/(n+2)}_x(I \times \R^n)}\\
&+ \| |\nabla u_\lo| |u_\lo|^{4/n} \|_{L^2_t L^{2n/(n+2)}_x(I \times \R^n)}
\lesssim_{m} \delta^\eps.
\end{align*}
By \eqref{mass-bound}, Corollary \ref{basic}, and Proposition \ref{highfreq}, we estimate
\begin{align*}
\| |u_\hi| |u_\lo|^{4/n} \|_{L^2_t L^{2n/(n+2)}_x(I \times \R^n)}
&\lesssim \|u_\lo\|_{L_t^\infty L_x^2(I \times \R^n)}^{4/n} \|u_\hi\|_{\S(I \times \R^n)}
\lesssim_{m}\delta + \delta^{4/n} \lesssim_{m} \delta^\eps\\
\| |u_\hi| |u_\hi|^{4/n} \|_{L^2_t L^{2n/(n+2)}_x(I \times \R^n)}
&\lesssim \|u_\hi\|_{L_t^\infty L_x^2(I \times \R^n)}^{4/n} \|u_\hi\|_{\S(I \times \R^n)}
\lesssim_{m}\delta + \delta^{4/n} \lesssim_{m} \delta^\eps.
\end{align*}
Hence, it remains to prove
\begin{equation}\label{uu}
\| |\nabla u_\lo| |u_\lo|^{4/n} \|_{L^2_t L^{2n/(n+2)}_x(I \times \R^n)}
\lesssim_{m} \delta^\eps.
\end{equation}
From \eqref{boot} and \eqref{qidef}, we have
$$ \| |x|^{-(1+\eps)/2} \nabla u_\lo \|_{L^2_{t,x}(I \times \R^n)} \lesssim \eta^{1/2} \lesssim 1$$
and thus, by radial Sobolev embedding (Corollary~\ref{radial embedding}),
$$ \| |\nabla|^{-\frac{1+\eps}{2(n-1)}} \nabla u_\lo \|_{L^2_t L^{q}_x(I \times \R^n)}
\lesssim\| |x|^{-(1+\eps)/2} \nabla u_\lo \|_{L^2_{t,x}(I \times \R^n)}
\lesssim 1
$$
where $q := 2(n-1)/(n-2-\eps)$.
Applying Bernstein we conclude that
$$ \| \nabla u_\lo \|_{L^2_t L^q_x(I \times \R^n)}  \lesssim 1.$$
As $q < 2n/(n-2)$, by H\"older we get
$$
\| |\nabla u_\lo| |u_\lo|^{4/n} \|_{L^2_t L^{2n/(n+2)}_x(I \times \R^n)}
\lesssim \| u_\lo \|_{L^\infty_t L^p_x(I \times \R^n)}^{4/n}
$$
for some $p > 2$.  But from \eqref{mass-bound}, \eqref{delta-bound}, Sobolev embedding, and Bernstein, we have
$$ \| u_\lo \|_{L^\infty_t L^p_x(I \times \R^n)} \lesssim_m \delta^c$$
for some $c > 0$, and so the contribution of this term is acceptable.

To complete the proof of Proposition~\ref{mor}, it remains to show that in dimensions $n=3,4$ we have
$$
\int_I \int_{\R^n} |G(t,x)| \frac{|u_\lo(t,x)|}{\langle x \rangle} \ dx dt \lesssim_{m} \delta^\eps,
$$
or equivalently (see the decomposition of $G$),
\begin{align}\label{ts-sd}
\|\langle x \rangle^{-1}|u_\lo| P_\lo O(|u_\hi||u_\lo|^{4/n})\|_{L_{t,x}^1}+\|\langle x \rangle^{-1}|u_\lo| P_\lo O(|u_\hi||u_\hi|^{4/n})\|_{L_{t,x}^1}\notag\\
+\|\langle x \rangle^{-1}|u_\lo|\Delta^{-1} \nabla \cdot P_\hi O( |u_\lo|^{4/n} |\nabla u_\lo| )_{L_{t,x}^1}
\lesssim_{m} \delta^\eps,
\end{align}
where all spacetime norms are on $I \times \R^n$.

To estimate the first term on the left-hand side of \eqref{ts-sd}, we use H\"older, Definition~\ref{sdef}, \eqref{lo-sd}, Proposition~\ref{highfreq},
and the fact that $|\nabla|^{-\frac{1-\eps}2}P_{>1/100}$ is bounded on $\S$ and $P_\lo$ is bounded on $L^2_x(|x|^{-\frac{n(1+\eps)+8}{n+4}})$:
\begin{align*}
\|\langle x \rangle^{-1}|u_\lo| P_\lo O(|u_\hi||u_\lo|^{4/n})\|_{L_{t,x}^1}
&\lesssim \||x|^{-\frac{n(1-\eps)}{2(n+4)}}u_\lo\|_{L_{t,x}^{2(n+4)/n}}^{(n+4)/n} \||x|^{-\frac{1+\eps}2}u_\hi\|_{L_{t,x}^2}  \\
&\lesssim \delta^{4/n} \||\nabla|^{-\frac{1-\eps}2}u_\hi\|_{\S}\\
&\lesssim \delta^{4/n} \|u_\hi\|_{\S}\\
&\lesssim \delta^{4/n}(\delta+\delta^{4/n}).
\end{align*}
Similarly, using \eqref{lo-so} instead of Proposition~\ref{highfreq} and the fact that $\Delta^{-1} \nabla \cdot P_\hi$ is bounded on
$L^2_x(|x|^{-\frac{n(1+\eps)+8}{n+4}})$, we estimate the third term on the left-hand side of \eqref{ts-sd} as follows:
\begin{align*}
\|\langle x \rangle^{-1}|u_\lo|\Delta^{-1} \nabla \cdot P_\hi O( |u_\lo|^{4/n} & |\nabla u_\lo| )_{L_{t,x}^1}\\
&\lesssim \||x|^{-\frac{n(1-\eps)}{2(n+4)}}u_\lo\|_{L_{t,x}^{2(n+4)/n}}^{(n+4)/n} \||x|^{-\frac{1+\eps}2}\nabla u_\lo\|_{L_{t,x}^2}  \\
&\lesssim \delta^{4/n}\||\nabla|^{-\frac{1-\eps}2}\nabla u_\lo\|_{\S}\\
&\lesssim \delta^{4/n}.
\end{align*}
To estimate the second term on the left-hand side on \eqref{ts-sd}, we use H\"older, Corollary~\ref{radial embedding} (or Hardy's inequality), Bernstein,
\eqref{delta-bound}, and the fact that $(\frac{n+4}n, \frac{2(n+4)}n)$ is a Schr\"odinger admissible pair in dimensions $n=3,4$, as well as the fact that
$P_\lo$ is bounded on $L_x^{2(n+4)/n}$ to get
\begin{align*}
\|\langle x \rangle^{-1}|u_\lo| P_\lo O(|u_\hi||u_\hi|^{4/n})\|_{L_{t,x}^1}
&\lesssim \||x|^{-1}u_\lo\|_{L_t^\infty L_x^2} \|u_\hi\|_{L_t^{(n+4)/n} L_x^{2(n+4)/n}}^{(n+4)/n}\\
&\lesssim \|\nabla u_\lo\|_{L_t^\infty L_x^2} \|u_\hi\|_{\S}^{(n+4)/n}\\
&\lesssim \delta(\delta+\delta^{4/n})^{(n+4)/n}.
\end{align*}

Putting everything together, we derive \eqref{ts-sd}.  The proof of Proposition \ref{mor} is now complete.

\subsection{Proof of Proposition \ref{decay}}\label{decaysec}

With Proposition \ref{mor} and all the above tools it is now an easy matter to establish
Proposition \ref{decay}.  Let $u$ be as in Proposition~\ref{decay}, and let $\eta > 0$ be an arbitrary
small quantity. From mass conservation we have \eqref{mass-bound} for some $m$. From Proposition
\ref{mor3} and a continuity argument we know that if \eqref{delta-bound} holds for some
sufficiently small $\delta$ (depending on $\eta$), then \eqref{boot} holds; applying Propositions
\ref{lowfreq} and \ref{highfreq} we then conclude the estimates
\begin{align*}
\| P_{\geq 1} u \|_{\S(I \times \R^n)} &\lesssim (\delta+\delta^{4/n}) \\
\| |\nabla|^{-(1-\eps)/2} \nabla P_{< 1} u \|_{\S(I \times \R^n)} &\lesssim \eta^{1/2}.
\end{align*}
One can now rescale these statements using \eqref{scaling}, replacing the role of the frequency $1$
by any other frequency $N$ (replacing $\nabla$ with $\nabla/N$).  Applying Lemma \ref{massdecay} to
make $\delta$ and $\eta$ arbitrarily small as $N \to \infty$ we then obtain Proposition
\ref{decay}.

\section{Proof of mass non-evacuation}\label{mass-sec}

We now prove Proposition \ref{nonevac}.  Assume for contradiction that we have a solution $u: I
\times \R^n \to \C$ obeying the hypotheses of that proposition, but such that
\begin{equation}\label{cascade}
\inf_{t \in I} N(t) = 0.
\end{equation}
Informally, \eqref{cascade} means that we have an unbounded cascade of mass from high frequencies
to low frequencies.  It turns out that this cascade, combined with perturbation theory and the Morawetz
estimate, allows us to improve the qualitative decay in Proposition \ref{decay} substantially:

\begin{proposition}[Cascade implies regularity]\label{cascadia}
Let the hypotheses be as above.  Then
\begin{equation}\label{cas}
\limsup_{N \to \infty} N^{(3-\eps)/2} \| u_{\geq N} \|_{\S} <\infty.
\end{equation}
Here and in the rest of this section, all spacetime norms are understood to be on the domain $I \times \R^n$.
\end{proposition}

\begin{proof}
From \eqref{massconserv} and Proposition \ref{decay}, we know that there exists $0 < m< +\infty$ such that
\begin{equation}\label{mass-u}
\| u \|_{L^\infty_t L^2_x} \leq m.
\end{equation}
Now let $\eta$ be a small number to be chosen later.  Then, by Proposition~\ref{decay} there exists $N_* > 0$ such that
$$\| u_{\geq N_*} \|_{\S}+ \frac{1}{N_*^{(1+\eps)/2}} \| |\nabla|^{-(1-\eps)/2} \nabla u_{< N_*} \|_{\S} \leq \eta.$$
Applying the scaling \eqref{scaling} (which does not affect qualitative hypotheses such as \eqref{freqbound} or \eqref{cascade},
the mass bound \eqref{mass-u}, or the qualitative conclusion \eqref{cas}), we may take $N_* = 1$; thus
\begin{align}
\| |\nabla|^{-(1-\eps)/2} \nabla u_{<1} \|_{\S} &\leq \eta \label{lo-u2}\\
\| u_{\geq 1} \|_{\S} &\leq \eta.\label{hi-u}
\end{align}

For any $\delta > 0$, let $P(\delta)$ denote the assertion that
$$ \| u_{\geq N} \|_{\S} \leq \eta N^{-(3-\eps)/2} + \delta, \, \hbox{ for all } N \geq 1.$$
By \eqref{hi-u}, we see that $P(\delta)$ is true for $\delta$ equal to $\eta$.  We now claim that if $\eta$ is sufficiently small,
then we have the bootstrap implication
\begin{equation}\label{boot-delta}
P(\delta) \implies P(\delta/2), \, \hbox{ for all } 0 < \delta \leq \eta.
\end{equation}
Iterating this to send $\delta$ to zero, we conclude that
$$ \| u_{\geq N} \|_{\S} \leq \eta N^{-(3-\eps)/2}, \, \hbox{ for all } N \geq 1$$
and the claim follows.

It remains to prove \eqref{boot-delta}.  Let $0 < \delta \leq \eta$ be such that $P(\delta)$ holds.
Let $N_0\geq 1$ be a dyadic integer such that $\eta N_0^{-(3-\eps)/2} \sim \delta$; then,
\begin{equation}\label{oon}
 \| u_{\geq N} \|_{\S} \lesssim \eta N^{-(3-\eps)/2} \quad \text{for all} \quad 1 \leq N \leq N_0
\end{equation}
and
\begin{align}\label{u>N_0}
\| u_{\geq N_0} \|_{\S} \lesssim \delta.
\end{align}

Now let $N\geq 1$ be arbitrary.  Applying $P_{\geq N}$ to \eqref{nls} we have
$$ (i\partial_t + \Delta) u_{\geq N} = P_{\geq N} F( u )$$
and hence, by Proposition \ref{strich}, we estimate
$$ \| u_{\geq N} \|_{\S} \lesssim \| u_{\geq N}(t_0) \|_{L^2_x(\R^n)} + \| P_{\geq N} F( u ) \|_{\N}$$
for any $t_0 \in I$.  From \eqref{cascade} and Definition \ref{apdef} we see that
$$ \inf_{t_0 \in I} \| u_{\geq N}(t_0) \|_{L^2_x(\R^n)} = 0.$$
Thus,
$$ \| u_{\geq N} \|_{\S} \lesssim \| P_{\geq N} F( u ) \|_{\N}.$$
We split
$$
 F(u) = F(u_{<N_0}) + O( |u_{\geq N_0}| |u_{<1}|^{4/n} ) + O( |u_{\geq N_0}| |u_{\geq 1}|^{4/n} )$$
and write
\begin{align*}
F(u_{<N_0}) = \Delta^{-1} \nabla \cdot \nabla F(u_{<N_0}) = \Delta^{-1} \nabla \cdot O( |u_{<N_0}|^{4/n} |\nabla u_{<N_0}| ).
\end{align*}
We conclude
\begin{align}
\| u_{\geq N} \|_{\S} &\lesssim
\| \Delta^{-1} \nabla \cdot P_{\geq N} O( |u_{<N_0}|^{4/n} |\nabla u_{<N_0}| )\|_{\N} \label{h1}\\
&\quad + \| P_{\geq N} O(|u_{\geq N_0}| |u_{<1}|^{4/n} )  \|_{\N}\label{h3}\\
&\quad + \| P_{\geq N} O(|u_{\geq N_0}| |u_{\geq 1}|^{4/n} ) \|_{\N}\label{h4}.
\end{align}

To estimate \eqref{h1}, we discard the H\"ormander-Mikhlin multiplier
$N^{(3-\eps)/2} \Delta^{-1} \nabla \cdot P_{\geq N} |\nabla|^{-(1-\eps)/2}$ and use the second estimate in Proposition \ref{refined}:
\begin{align*}
\| \Delta^{-1} \nabla \cdot P_{\geq N} & O( |u_{<N_0}|^{4/n} |\nabla u_{<N_0}| )\|_{\N}\\
&\lesssim N^{-(3-\eps)/2} \| |\nabla|^{(1-\eps)/2} O( |u_{<N_0}|^{4/n} |\nabla u_{<N_0}| ) \|_{\N} \\
&\lesssim N^{-(3-\eps)/2} \| |\nabla|^{\frac{3}{4}(1-\eps)} u_{<N_0} \|_{L^\infty_t L^2_x}^{4/n}
\| |\nabla|^{(1-\eps)(\frac{1}{2} - \frac{3}{n})}  \nabla u_{<N_0} \|_{\S}.
\end{align*}
From \eqref{lo-u2}, \eqref{oon}, and Bernstein, we see that
\begin{align*}
\| |\nabla|^{\frac{3}{4}(1-\eps)} u_{<N_0} \|_{L^\infty_t L^2_x}
&\lesssim \| |\nabla|^{\frac{3}{4}(1-\eps)} u_{<1} \|_{\S} + \sum_{1\leq M<N_0}\| |\nabla|^{\frac{3}{4}(1-\eps)} u_{M} \|_{\S}\\
&\lesssim \eta + \eta\sum_{1\leq M<N_0}M^{\frac{3}{4}(1-\eps)}M^{-\frac{3-\eps}2}\\
&\lesssim \eta.
\end{align*}
Similarly,
\begin{align*}
\| |\nabla|^{(1-\eps)(\frac{1}{2} - \frac{3}{n})}  \nabla u_{<N_0} \|_{\S}
&\lesssim \| |\nabla|^{(1-\eps)(\frac{1}{2} - \frac{3}{n})}  \nabla u_{<1} \|_{\S}+ \sum_{1\leq M<N_0}\| |\nabla|^{(1-\eps)(\frac{1}{2} - \frac{3}{n})}  \nabla u_M \|_{\S}\\
&\lesssim \eta + \eta\sum_{1\leq M<N_0}M^{1+(1-\eps)(\frac{1}{2} - \frac{3}{n})}M^{-\frac{3-\eps}2}\\
&\lesssim \eta.
\end{align*}
We thus conclude
$$ \eqref{h1} \lesssim \eta^{4/n} \eta N^{-(3-\eps)/2}. $$

To estimate \eqref{h3}, we discard the H\"ormander-Mikhlin multiplier $N^{\frac{1-\eps}2} P_{\geq N} |\nabla|^{-\frac{1-\eps}2}$ and use the first
estimate in Proposition~\ref{refined}:
\begin{align*}
\| P_{\geq N} O(|u_{\geq N_0}| |u_{<1}|^{4/n} ) \|_{\N}
&\lesssim N^{-(1-\eps)/2}\| |\nabla|^{(1-\eps)/2} O(|u_{\geq N_0}| |u_{<1}|^{4/n} ) \|_{\N} \\
&\lesssim N^{-(1-\eps)/2} \| |\nabla|^{\frac{n}{4}(1-\eps)} u_{<1} \|_{L^\infty_t L^2_x}^{4/n} \| |\nabla|^{-(1-\eps)/2} u_{\geq N_0} \|_{\S}.
\end{align*}
From the boundedness of the H\"ormander-Mikhlin multiplier $N_0^{\frac{1-\eps}2} P_{\geq N_0} |\nabla|^{-\frac{1-\eps}2}$ on $\S$,
and \eqref{u>N_0}, we get
$$ \| |\nabla|^{-(1-\eps)/2} u_{\geq N_0} \|_{\S} \lesssim N_0^{-(1-\eps)/2} \delta,$$
while from \eqref{lo-u2} and Bernstein, we get
\begin{align*}
\| |\nabla|^{\frac{n}{4}(1-\eps)} u_{<1} \|_{L^\infty_t L^2_x}\lesssim \eta.
\end{align*}
Putting all these together we obtain
$$ \eqref{h3} \lesssim N^{-(1-\eps)/2}N_0^{-(1-\eps)/2}\eta^{4/n} \delta\lesssim \eta^{4/n} \delta.$$

Finally, to estimate \eqref{h4}, we discard the H\"ormander-Mikhlin multiplier $P_{\geq N}$ and use
Corollary \ref{basic} followed by \eqref{hi-u} and \eqref{u>N_0}:
\begin{align*}
\| P_{\geq N}O(|u_{\geq N_0}| |u_{\geq 1}|^{4/n} ) \|_{\N}
&\lesssim \| O(|u_{\geq N_0}| |u_{\geq 1}|^{4/n} )\|_{\N} \\
&\lesssim \| u_{\geq 1} \|_{L^\infty_t L^2_x}^{4/n}\| u_{\geq N_0} \|_{\S} \\
&\lesssim \eta^{4/n} \delta .
\end{align*}

Combining all these estimates we conclude that
$$
\| u_{\geq N} \|_{\S} \lesssim \eta^{4/n}(\eta N^{-(3-\eps)/2} + \delta)
$$
for all $N\geq 1$, which (for $\eta$ small) implies $P(\delta/2)$ as desired.

This concludes the proof of Proposition \ref{cascadia}.
\end{proof}

We can now combine the regularity given by \eqref{cas} with the high-to-low frequency cascade \eqref{cascade} to contradict energy conservation.
From Proposition~\ref{cascadia} we have
$$ \| P_N u \|_{L^\infty_t L^2_x} \lesssim_u N^{-(3-\eps)/2}$$
for all sufficiently large $N$.  Of course, from mass conservation we also have
$$ \| P_N u \|_{L^\infty_t L^2_x} \lesssim_u 1$$
for all $N$.  So, in particular, $u \in L^\infty_t H^1_x$ and
$$ \| P_N \nabla u \|_{L^\infty_t L^2_x}^2 \lesssim_u \min( N^2, N^{-1+\eps} ).$$
Now, from \eqref{cascade} we can find a sequence $t_i \in I$ with $N(t_i) \to 0$.  Then, by Definition~\ref{apdef} we see that
$$ \lim_{i \to \infty} \| P_N \nabla u(t_i) \|_{L^\infty_t L^2_x}^2 = 0$$
for all $N$.  On the other hand, we have just established that
$$ \| P_N \nabla u(t_i) \|_{L^\infty_t L^2_x}^2 \lesssim_u \min( N^2, N^{-1+\eps} ).$$
Observe that $\min( N^2, N^{-1+\eps} )$ is absolutely summable as $N$ ranges over dyadic frequencies.  Thus, by the dominated convergence theorem we have
$$ \lim_{i \to \infty} \sum_N \| P_N \nabla u(t_i) \|_{L^\infty_t L^2_x}^2 = 0$$
and hence by orthogonality,
$$ \lim_{i \to \infty} \| \nabla u(t_i) \|_{L^\infty_t L^2_x} = 0$$
From mass conservation and the Gagliardo-Nirenberg inequality, we then get
$$ \lim_{i \to \infty} \| u(t_i) \|_{L^\infty_t L^{2(n+2)/n}_x} = 0$$
and hence,
$$ \lim_{i \to \infty} E(u(t_i)) = 0$$
where $E$ is the \emph{energy}
$$ E(u(t_i)) := \int_{\R^n} \frac{1}{2} |\nabla u(t_i,x)|^2 + \frac{n}{2(n+2)} |u(t,x)|^{2(n+2)/n}\ dx.$$
But $H^1$ solutions to \eqref{nls} conserve the energy (see e.g. \cite{caz}), and hence $E(u(t)) = 0$
for all $t \in I$.  In other words, $u$ is identically zero, contradicting the hypotheses of
Proposition \ref{nonevac}.  This shows that \eqref{cascade} cannot occur and
Proposition~\ref{nonevac} follows.

\appendix

\section{Weighted Sobolev estimates}\label{appendix}

In this appendix we collect some estimates of Sobolev type which will allow us to manipulate power
weights such as $|x|^{\pm \alpha}$ or $\langle x \rangle^{\pm \alpha}$ relatively easily.  We begin
with a bilinear form estimate.

\begin{lemma}[Bilinear form estimate]\label{bil}
Let $n \geq 1$, $1 \leq p,q \leq \infty$ be such that $\frac{1}{p} + \frac{1}{q} \geq 1$, and
let $f \in L^p(\R^n)$ and $g \in L^q(\R^n)$.  Let $\alpha, \beta \in \R$ obey the scaling condition
$$ \alpha + \beta = -\frac{n}{p'} - \frac{n}{q'}.$$
\begin{itemize}
\item If $\alpha > -\frac{n}{p'}$ (thus $\beta < -\frac{n}{q'}$), then
$$ \int\int_{|x| \lesssim |y|} |x|^\alpha |y|^\beta |f(x)| |g(y)| \ dx dy \lesssim_{\alpha,\beta,p,q} \|f\|_{L^p(\R^n)} \|g\|_{L^q(\R^n)}.$$
\item If $\alpha < -\frac{n}{p'}$ (thus $\beta > -\frac{n}{q'}$), then
$$ \int\int_{|y| \lesssim |x|} |x|^\alpha |y|^\beta |f(x)| |g(y)| \ dx dy \lesssim_{\alpha,\beta,p,q} \|f\|_{L^p(\R^n)} \|g\|_{L^q(\R^n)}.$$
\end{itemize}
Here $p',q'$ are the dual exponents to $p,q$, that is, $\frac{1}{p} + \frac{1}{p'} = \frac{1}{q} + \frac{1}{q'} = 1$.
\end{lemma}

\begin{proof} From H\"older's inequality, we see that
\begin{align*}
\int_{|y| \sim R} \int_{|x| \sim 2^m R} & |x|^\alpha |y|^\beta f(x) g(y) \ dx dy\\
&\lesssim_{\alpha,\beta,p,q} (2^m R)^\alpha R^\beta (2^m R)^{n/p'} R^{n/q'}
\|f\|_{L^p(|x| \sim 2^m R)} \|g\|_{L^q(|y| \sim R)}
\end{align*}
for any $R > 0$ and $m \in \Z$.  The scaling
condition ensures that the powers of $R$ cancel. Summing over dyadic $R$ and using H\"older (taking
advantage of the hypothesis $\frac{1}{p} + \frac{1}{q} \geq 1$), we conclude that
$$ \int\int_{|x| \sim 2^m |y|} |x|^\alpha |y|^\beta f(x) g(y) \ dx dy \lesssim_{\alpha,\beta,p,q} 2^{m(\alpha + n/p')}
\|f\|_{L^p(\R^n)} \|g\|_{L^q(\R^n)}.$$
Summing over either negative or positive $m$ then yields the claim.
\end{proof}

Now, we insert a fractional integral weight $\frac{1}{|x-y|^{n-s}}$ and assume spherical symmetry.

\begin{lemma}[Hardy-Littlewood-Sobolev estimate]\label{hs}
Let $1 \leq p, q \leq \infty$, $n \geq 1$, $0 < s < n$, and $\alpha,\beta \in \R$ obey the
conditions
\begin{align*}
\alpha &>  - \frac{n}{p'} \\
\beta &> -\frac{n}{q'}\\
1 \leq \frac{1}{p} + \frac{1}{q} &\leq 1+s
\end{align*}
and the scaling condition
$$ \alpha + \beta - n + s = - \frac{n}{p'} - \frac{n}{q'},$$
with at most one of the equalities
$$ p=1, p=\infty, q=1, q=\infty, \frac{1}{p} + \frac{1}{q} = 1 + s$$
holding. Let $f \in L^p(\R^n)$ and $g \in L^q(\R^n)$ be spherically
symmetric.  Then
$$ \int_{\R^n} \int_{\R^n} \frac{|x|^\alpha |y|^\beta}{|x-y|^{n-s}} |f(x)| |g(y)| \ dx dy \lesssim_{\alpha,\beta,p,q} \|f\|_{L^p(\R^n)} \|g\|_{L^q(\R^n)}.$$
\end{lemma}

\begin{proof}
From Lemma \ref{bil} we can already dispose of the portion of the integral where $|x| \leq |y|/2$
or $|y| \leq |x|/2$, since in those cases one can replace $\frac{1}{|x-y|^{n-s}}$ by
$\frac{1}{|y|^{n-s}}$ or $\frac{1}{|x|^{n-s}}$, respectively.  Thus, it suffices to show that
$$ \int\int_{|x| \sim |y|} \frac{|x|^\alpha |y|^\beta}{|x-y|^{n-s}} |f(x)| |g(y)| \ dx dy \lesssim_{\alpha,\beta,p,q} \|f\|_{L^p(\R^n)} \|g\|_{L^q(\R^n)}.$$
Since $\frac{1}{p}+\frac{1}{q} \geq 1$, we see that it will suffice to show that
$$ \int_{|y| \sim R} \int_{|x| \sim R} \frac{|x|^\alpha |y|^\beta}{|x-y|^{n-s}} |f(x)| |g(y)| \ dx dy \lesssim_{\alpha,\beta,p,q} \|f\|_{L^p(|x| \sim R)} \|g\|_{L^q(|y| \sim R)}$$
for all $R > 0$, since one can then sum over dyadic $R$ and use H\"older.  The scaling condition
allows us to normalize $R=1$; thus, we reduce to showing
$$ \int_{|y| \sim 1} \int_{|x| \sim 1} \frac{1}{|x-y|^{n-s}} |f(x)| |g(y)| \ dx dy \lesssim_{p,q} \|f\|_{L^p(|x| \sim 1)} \|g\|_{L^q(|y| \sim 1)}.$$
Suppose first that $s < 1$.  Then by switching to polar co-ordinates $f(x) = F(|x|)$, $g(y) = G(|y|)$ we reduce to the one-dimensional estimate
$$ \int_{|y| \sim 1} \int_{|x| \sim 1} \frac{1}{|x-y|^{1-s}} |F(x)| |G(y)| \ dx dy \lesssim_{p,q} \|f\|_{L^p(|x| \sim 1)} \|g\|_{L^q(|y| \sim 1)}$$
which follows from the classical Hardy-Littlewood-Sobolev inequality.  If instead $s \geq 1$, then we (at worst) reduce to
$$ \int_{|y| \sim 1} \int_{|x| \sim 1} (\log \frac{1}{|x-y|}) |F(x)| |G(y)| \ dx dy \lesssim_{p,q} \|f\|_{L^p(|x| \sim 1)} \|g\|_{L^q(|y| \sim 1)}$$
which follows from Young's inequality.
\end{proof}

By using duality we now obtain

\begin{corollary}[Radial Sobolev inequality]\label{radial embedding}
Let $n, \alpha, \beta, p, q, s$ be as in Lemma \ref{hs}.  Then for any spherically symmetric $u: \R^n \to \C$, we have
$$ \| |x|^\beta u \|_{L^{q'}(\R^n)} \lesssim_{\alpha,\beta,p,q,s}
\| |x|^{-\alpha} |\nabla|^s u \|_{L^p(\R^n)}.$$
\end{corollary}

\begin{proof} Write $f := |x|^{-\alpha} |\nabla|^s u$; we are trying to show
$$ \| |x|^\beta |\nabla|^{-s}(|x|^\alpha f) \|_{L^{q'}(\R^n)} \lesssim_{\alpha,\beta,p,q,s} \| f \|_{L^p(\R^n)}$$
which by duality is equivalent to
$$ \int_{\R^n} g(x) |x|^\beta |\nabla|^{-s}(|x|^\alpha f)\ dx
\lesssim_{\alpha,\beta,p,q,s} \|f\|_{L^p(\R^n)} \|g\|_{L^q(\R^n)}.$$
But this follows from Lemma \ref{hs}.
\end{proof}

Finally, we will need to switch between homogeneous power weights $|x|^{-\alpha}$ and inhomogeneous
power weights $\langle x \rangle^{-\alpha}$.  In one direction this is trivial since $\langle x
\rangle^{-\alpha} \lesssim_\alpha |x|^{-\alpha}$ for $\alpha \geq 0$.  If the frequency of the
function is localized to $O(N)$, one expects to be able to reverse this inequality up to a factor
of $\langle N \rangle^\alpha$ due to the spatial uncertainty of $O(1/N)$.  More precisely, we have

\begin{lemma}[Uncertainty principle]\label{smoother} If $f: \R^n \to \C$, $1 < p < \infty$,
$0 < \alpha < n/p$, and $N > 0$, then
$$ \| |x|^{-\alpha} P_{<N} f \|_{L^p(\R^n)} \lesssim_{\alpha,p} \langle N \rangle^\alpha \| \langle x \rangle^{-\alpha} f \|_{L^p(\R^n)}.$$
\end{lemma}

\begin{proof}
We may assume that $N \geq 1$, since the case for smaller $N$ follows from the $N=2$ case by using
the factorization $P_{<N} = P_{<N} P_{<2}$ and then discarding $P_{<N}$ using Lemma \ref{calculus}.

On the other hand, for $N \geq 1$, the claim would follow immediately from
\begin{align}\label{smoother-reduction}
\| |Nx|^{-\alpha} P_{<N} f \|_{L^p(\R^n)} \lesssim \| \langle Nx \rangle^{-\alpha} f \|_{L^p(\R^n)},
\end{align}
as $\langle x\rangle^{-\alpha}$ can be estimated from below by $\langle Nx\rangle^{-\alpha}$.

It thus suffices to prove \eqref{smoother-reduction}; rescaling by $N$, we may normalize $N=1$.  If $f$ is supported on the set $\{ x: |x| \geq 1\}$,
then we estimate $\langle  x \rangle^{-\alpha}$ from below by $|x|^{-\alpha}$, and \eqref{smoother-reduction} follows from Lemma \ref{calculus}
since $P_{<1}$ is a H\"ormander-Mikhlin  multiplier.  So, we may reduce to the case when $f$ is
supported on the ball $\{ x: |x| \leq 1 \}$.  But then, a direct computation using the convolution
kernel of $P_{<1}$ and H\"older, establishes the pointwise estimate
$$ P_{<1} f(x) = O( \langle x \rangle^{-100n} \|f\|_{L^p(\R^n)} )$$
and \eqref{smoother-reduction} follows.
\end{proof}


\begin{thebibliography}{10}

\bibitem{begout}
P. Begout, A. Vargas, \emph{Mass concentration phenomena for the $L^2$-critical nonlinear Schr\"odinger equation}, preprint.

\bibitem{bourg.2d}
J. Bourgain, \emph{Refinements of Strichartz Inequality and Applications
to 2d-NLS With Critical Nonlinearity}, Inter. Math. Res. Not.  (1998), 253--284.

\bibitem{borg:scatter}
J. Bourgain, \emph{Global well-posedness of defocusing 3D critical NLS in the radial case}, JAMS
\textbf{12} (1999), 145--171.

\bibitem{borg:book}
J. Bourgain, \emph{New global well-posedness results for non-linear Schr\"odinger equations}, AMS
Publications (1999).

\bibitem{ck}
R. Carles, S. Keraani, \emph{On the role of quadratic oscillations in nonlinear Schr\"odinger equation II.  The $L^2$-critical case}, preprint.

\bibitem{cwI}
T. Cazenave, F.B. Weissler, \emph{Critical nonlinear Schr\"odinger Equation}, Non. Anal. TMA \textbf{14}
(1990), 807--836.

\bibitem{caz}
T. Cazenave, \textit{Semilinear Schr\"odinger equations,} Courant
Lecture Notes in Mathematics, \textbf{10}, American Mathematical
Society, 2003.

\bibitem{cct:ill} M. Christ, J. Colliander, T. Tao, \emph{Ill-posedness for nonlinear Schr\"odinger and wave equations}, preprint
\texttt{math.AP/0311048}.

\bibitem{ckstt:gwp}
J. Colliander, M. Keel, G. Staffilani, H. Takaoka, T. Tao, \emph{Global well-posedness and scattering for the energy-critical
nonlinear Schr\"odinger equation in $\R^3$}, to appear Annals of Math.

\bibitem{ckstt:french}
J. Colliander, M. Keel, G. Staffilani, H. Takaoka, T. Tao, \emph{Existence globale et diffusion pour
l'\'equation de
  Schr\"odinger nonlin\'eaire r\'epulsive cubique sur $\R^3$ en
  dessous l'espace d'\'energie},
Journ\'ees ``\'Equations aux D\'eriv\'ees Partielles''
              (Forges-les-Eaux, 2002), Exp. No. X, 14, 2002.

\bibitem{glassey}
R. T. Glassey, \emph{On the blowing up of solution to the Cauchy problem for nonlinear Schr\"odinger
operators}, J. Math. Phys. \textbf{8} (1977), 1794--1797.

\bibitem{gv:scatter}
J. Ginibre, G. Velo, \emph{Scattering theory in the energy space for a class of nonlinear Schr\"odinger
equations}, J. Math. Pure. Appl. \textbf{64} (1985), 363--401.

\bibitem{hayashi:tsutsumi}
N. Hayashi, Y. Tsutsumi, \emph{Remarks on the scattering problem for nonlinear Schr\"odinger equations},
Diff. Eq. Math. Phys. (1986), 162--168.

\bibitem{kato}
T. Kato, \emph{On nonlinear Schr\"odinger equations},  Ann. Inst. H. Poincar\'e Phys. Theor. \textbf{46} (1987), 113--129.

\bibitem{katounique}
T. Kato, \emph{On nonlinear Schr\"odinger equations, II.  $H^s$-solutions and unconditional well-posedness},
J. d'Analyse. Math. \textbf{67} (1995), 281--306.

\bibitem{tao:keel}
M. Keel, T. Tao, \emph{Endpoint Strichartz Estimates}, Amer. Math. J. \textbf{120} (1998), 955--980.

\bibitem{merlekenig}
C. Kenig, F. Merle, \emph{Global well-posedness, scattering, and blowup for the energy-critical, focusing, non-linear Schr\"odinger equation in the radial case}, preprint.

\bibitem{keraani}
S. Keraani, \emph{On the blow-up phenomenon of the critical nonlinear Schr\"odinger equation}, preprint.

\bibitem{mt}
F. Merle, Y. Tsutsumi, \emph{$L^2$ concentration of blow-up solutions for the nonlinear Schr\"odinger equation with critical power nonlinearity},
J. Diff. Eq. \textbf{84} (1990), 205--214.

\bibitem{mv}
F. Merle, L. Vega, \emph{Compactness at blow-up time for $L^2$ solutions of the critical nonlinear Schr\"odinger equation in 2D},
Internat. Math. Res. Not. \textbf{8} (1998), 399--425.

\bibitem{merle}
F. Merle, \emph{Existence of blow-up solutions in the energy space for the critical generalized KdV equation}, J. Amer. Math. Soc.
\textbf{14} (2001), 555--578.

\bibitem{morawetz}
C. Morawetz, \emph{Time decay for the nonlinear Klein-Gordon equation}, Proc. Roy. Soc. A \textbf{306} (1968), 291--296.

\bibitem{RV}
E.~Ryckman, M.~Visan, \emph{Global well-posedness and scattering for the defocusing energy-critical
nonlinear Schr\"odinger equation in $\R^{1+4}$,} to appear Amer. J. Math.

\bibitem{stein:large}
E.~M. Stein, \emph{Harmonic Analysis}, Princeton University Press, 1993.

\bibitem{stein:weiss}
E. M. Stein, G. Weiss, \emph{Fractional integrals on n-dimensional Euclidean space}, J. Math. Mech. \textbf{7} (1958), 503--514.

\bibitem{SulSul}
C. Sulem, P.-L. Sulem, \emph{The nonlinear Schr\"odinger equation: Self-focusing and wave collapse}, Springer-Verlag, New York, (1999).

\bibitem{tao: gwp radial}
T. Tao, \emph{Global well-posedness and scattering for the higher-dimensional energy-critical non-linear
Schr\"odinger equation for radial data}, New York Journal of Math. \textbf{11} (2005), 57-80.

\bibitem{tao:bilinear}
T. Tao, \emph{A sharp bilinear restriction estimate for
paraboloids},  Geom. Funct. Anal. \textbf{13} (2003), 1359--1384.

\bibitem{tao-lens}
T. Tao, \emph{A pseudoconformal compactification of the nonlinear Schr\"odinger equation and applications}, preprin \texttt{math.AP/0606254}.

\bibitem{TV}
T. Tao, M.~Visan, \emph{Stability of energy-critical nonlinear Schr\"odinger equations in high
dimensions,} Electron. J. Diff. Eqns. \textbf{118} (2005), 1-28.

\bibitem{tvz} T. Tao, M.~Visan, X. Zhang, \emph{The nonlinear Schr\"odinger equation with combined power-type
nonlinearities,} preprint \texttt{math. AP/0511070}.

\bibitem{compact} T. Tao, M.~Visan, X. Zhang, \emph{Minimal-mass blowup solutions of the mass-critical NLS}, preprint \textbf{math.AP/0609690}.

\bibitem{tsutsumi}
Y. Tsutsumi, \emph{Scattering problem for nonlinear Schr\"odinger equations},
Ann. Inst. H. Poincar\'e Phys. Th\'eor. \textbf{43} (1985), 321--347.

\bibitem{tzirakis} N. Tzirakis, \emph{The Cauchy problem for the semilinear quintic Schr\"odinger equation in one dimension},
Differential Integral Equations \textbf{18} (2005), 947--960.

\bibitem{vilela}
M. Vilela, \emph{Regularity of solutions to the free Schr\"odinger equation with radial initial data}, Illinois J. Math. \textbf{45} (2001), 361--370.

\bibitem{thesis:art} M. Visan, \emph{The defocusing energy-critical nonlinear Schr\"odinger equation
in higher dimensions,} to appear Duke Math. J.

\bibitem{Monica:thesis} M. Visan, \emph{The defocusing energy-critical nonlinear Schr\"odinger equation
in dimensions five and higher,} Ph.D. Thesis, UCLA.

\end{thebibliography}
\end{document}